\documentclass[a4paper, 12pt]{article}

\title{\bf
 Mutation in triangulated categories and rigid Cohen-Macaulay modules
\\}
\author
{Osamu Iyama \\
{\small Nagoya University, }
{\small Nagoya 464-8602, Japan} \\
{\small \texttt{iyama@math.nagoya-u.ac.jp}} \\
\vspace{12pt} \\
Yuji Yoshino \\
{\small Okayama University, }
{\small Okayama 700-8530, Japan} \\
{\small \texttt{yoshino@math.okayama-u.ac.jp}}
}

\usepackage{amssymb}
\usepackage{amsmath}
\usepackage{amscd}
\usepackage{latexsym}
\usepackage{theorem}
\usepackage{color}

\theoremheaderfont{\scshape}
\newtheorem{thm}{\bfseries Theorem}[section]
\newtheorem{prop}[thm]{\bfseries Proposition}
\newtheorem{lemma}[thm]{\bfseries Lemma}
\newtheorem{cor}[thm]{\bfseries Corollary}

\theorembodyfont{\rmfamily}
\newtheorem{defn}[thm]{\bfseries Definition}
\newtheorem{ex}[thm]{Example}

\newtheorem{notation}[thm]{Notation}

\newtheorem{pf}{Proof.}

\def\Z{{\mathbb Z}}

\def\ZC{{\mathbb Z}/3{\mathbb Z}}

\def\gp{\mathfrak p}

\def\Hom{\mathrm{Hom}}
\def\pHom{\underline{\mathrm{Hom}}}
\def\pEnd{\underline{\mathrm{End}}}

\def\Ext{\mathrm{Ext}}

\def\End{\mathrm{End}}
\def\pEnd{\underline{\mathrm{End}}}

\def\Spec{\mathrm{Spec}}

\def\Kdim{\mathrm{dim}}

\def\rank{\mathrm{rank}}

\def\CM{\mathrm{CM}}

\def\CMR{\mathrm{CM}(R)}
\def\CML{\mathrm{CM}(\Lambda)}

\def\length{\mathrm{length}}
\def\pCMR{\underline{\mathrm{CM}}(R)}
\def\pCMA{\underline{\mathrm{CM}}(A)}
\def\pCML{\underline{\mathrm{CM}}(\Lambda)}

\def\GL{\mathrm{GL}}
\def\SL{\mathrm{SL}}

\def\qed{$\Box$}

\def\endm{\mathop{\rm End}\nolimits}

\def\pCMR{\underline{\mathrm{CM}}(R)}

\def\add{\mathop{\rm add }\nolimits}
\def\tt{\mathop{\cal T}\nolimits}
\def\sss{\mathbb{ S}}
\def\ppp{\mathbb{ P}}
\def\fff{\mathbb{ F}}
\def\ggg{\mathbb{ G}}
\def\hhh{\mathbb{ H}}

\def\cc{\mathop{\cal C}\nolimits}

\def\ind{\mathop{\rm ind}\nolimits}

\def\dd{\mathop{\cal D}\nolimits}
\def\uu{\mathop{\cal U}\nolimits}
\def\xx{\mathop{\cal X}\nolimits}
\def\yy{\mathop{\cal Y}\nolimits}
\def\zz{\mathop{\cal Z}\nolimits}
\def\Im{\mathop{\rm Im}\nolimits}
\def\Mod{\mathop{\rm Mod}\nolimits}

\def\height{\mathrm{ht}}

\def\mod{{\mathrm{mod }}}

\def\ang#1#2#3#4{\mathop{\stackrel{\begin{picture}(10,12)
\put(-6,4){\vector(1,1){10}}
\put(0,17){\scriptsize$#1$}
\put(8,14){\vector(1,-1){10}}
\put(-2,0){\scriptsize$#2$}
\put(-18,9){\tiny$#3$}
\put(15,9){\tiny$#4$}
\end{picture}}{\ \longrightarrow\ }}\nolimits}

\date{\empty}
\begin{document}

\maketitle


\begin{center}
{\it Dedicated to Professor Idun Reiten on the occasion of her 65th birthday}
\end{center}

\vskip1.5em
\begin{abstract}
We introduce the notion of mutation of $n$-cluster tilting subcategories in a triangulated category with Auslander-Reiten-Serre duality.
Using this idea, we are able to obtain the complete classifications of rigid Cohen-Macaulay modules over certain Veronese subrings.
\end{abstract}

\begin{center}
{\bf Contents}
\end{center}
\begin{itemize}
\item[1] Introduction
\item[2] Torsion theories and mutation pairs
\item[3] $n$-cluster tilting subcategories
\item[4] Subfactor triangulated categories
\item[5] Mutation of $n$-cluster tilting subcategories
\item[6] Subfactor abelian categories
\item[7] Application of Kac's theorem
\item[8] Gorenstein quotient singularities
\item[9] Proof of Theorems 1.2 and 1.3
\item[10] Non-commutative examples
\end{itemize}


\section{Introduction}

The theory of Auslander and Reiten is now a basic tool for the representation theory of finite dimensional algebras  and for the theory of Cohen-Macaulay modules over Cohen-Macaulay rings \cite{ARS,Yoshino1}. 
In particular, Auslander-Reiten theory enables us to give a complete description of   Cohen-Macaulay modules over a two dimensional simple singularity, which correspond to irreducible representations of the relevant finite group \cite{Auslander}.  
This is closely related to so-called McKay correspondence in algebraic geometry (cf. \cite{AV, GV, KV} etc.). For higher dimensional cases, many authors (cf. \cite{ItoNakajima, BKR, BK} etc.) have studied the resolution of singularities for the ring of invariants given by a finite group in $\SL(d,k)\  (d >2)$. 
However, from the viewpoint of module theory or representation theory, it is more natural to seek the classification of Cohen-Macaulay modules over the ring of invariants. 
However, a simple observation shows us that such rings of invariants often have wild representation type (see \ref{wild}). Hence, the task of providing a classification of all modules has been regarded as `hopeless' \cite{Drozd, Crawley-Boevey}.

In the present paper we focus our attention on the classification of rigid Cohen-Macaulay modules. 
We show that the method of mutation is one of the most promising methods for classifying them.
We give a complete classification of rigid Cohen-Macaulay modules in several special cases where the groups are given as 
 $G=\langle{\rm diag}(\omega,\omega,\omega)\rangle$ ($\omega$ is a primitive cubic root of unity) \cite{Yoshino2} and  $G = \langle{\rm diag}(-1, -1, -1, -1)\rangle$. 
Surprisingly, the general description can be given in terms of the root systems (see \ref{Kaccor}). 
We also note that rigid Cohen-Macaulay modules have deep relevance to the noncommutative crepant resolutions \cite{Iyama2, IyamaReiten} introduced by Van den Bergh \cite{VandenBergh1, VandenBergh2}.  

There are several basic theorems in Auslander-Reiten theory which play 
an essential role in representation theory. 
In \cite{Iyama1, Iyama2} the first author introduced the idea of maximal $(n-1)$-orthogonal subcategories for certain exact categories.
They were called $n$-cluster tilting subcategories in Keller and Reiten \cite{KR}, and we borrow this suggestive terminology in this paper.
Using this approach `higher dimensional' analogues of some results in Auslander-Reiten theory can be found: for example, the existence of higher Auslander-Reiten sequences. There is also an analogy for the one-one correspondence between Morita equivalence classes of finite dimensional algebras of finite representation type and classes of Auslander algebras. See \cite{Geiss, KR, EH, IyamaReiten} for further results.

In this paper, we study the corresponding idea not in abelian categories but in any triangulated categories \cite{Iyama3}. 
In particular, we introduce an idea that we call AR (=Auslander-Reiten) $(n+2)$-angles in $n$-cluster tilting subcategories of triangulated categories and we shall prove that they always exist. 
AR $(n+2)$-angles are a kind of complex which will take the place of Auslander-Reiten triangles \cite{Happel}. 
Moreover, strongly motivated by recent studies \cite{Buan, Geiss, KR} on cluster categories which we explain below, we introduce the idea of mutation. This is an operation on the set of $n$-cluster tilting subcategories, which is closely related to AR $(n+2)$-angles.
We believe this idea is an essential tool for the classification of rigid objects in a triangulated category. 

Now let us recall the history of mutations. 
Bern\v ste\u\i n, Gelfand and Ponomarev \cite{BGP} introduced reflection functors to the representation theory of quivers.
This was the starting point of later tilting theory \cite{APR, BB, Miyashita, Rickard, Keller}. 
The remarkable combinatorial aspect of tilting theory was first observed by Riedtmann and Schofield \cite{RS}.
Gorodentsev and Rudakov \cite{GR} (see also \cite{Rudakov}) made use of mutation when they classified the exceptional vector bundles on $\ppp^2$. 
Recently, Seidel and Thomas \cite{ST} constructed the twist autoequivalences on derived categories along similar lines. 
Generally speaking, mutations can be regarded as a categorical realization of Coxeter or braid groups. 

Recently, another strong motivation has come from the area of algebraic combinatorics.
Fomin and Zelevinsky \cite{FZ1, FZ2} introduced cluster algebras, which enjoy important combinatorial properties given in terms of the mutation for skew symmetric matrices.
As a categorification of cluster algebras, Buan, Marsh, Reineke, Reiten and Todorov \cite{Buan} and Caldero, Chapoton and Schiffler (for $A_n$ case) \cite{CCS}
introduced certain triangulated categories called cluster categories, which were defined as orbit categories \cite{Keller2} of the derived categories of quiver representations.
They introduced the mutation (called exchange there) of cluster tilting objects in these categories, which corresponds to the mutation of clusters in cluster algebras.
Geiss, Leclerc and Schr\"oer \cite{Geiss} applied mutation to the study of rigid modules over preprojective algebras and the coordinate rings of maximal unipotent subgroups of semisimple Lie groups. We refer the reader to \cite{BMR, BMR2, CK, KR, KZ, BM, Zhu, Thomas, Tabuada} for more recent developments in the study of cluster categories and their generalizations.
We believe that our theory will be useful for studying them as well as rigid Cohen-Macaulay modules. In fact, Keller and Reiten \cite{KR2} provide a specific relationship between them. Our results in this paper will be applied in \cite{BIRS} to study cluster structures for 2-Calabi-Yau triangulated categories. 


It is of the greatest interest to us that we are able to consider mutation as an aspect of the approximation theory initiated by the school of Auslander \cite{AS, AB, AR2.5}. 
The mutation of $n$-cluster tilting subcategories based on approximation theory that will be introduced in the present paper is a `higher dimensional' generalization of one given in \cite{Buan}  and  \cite{Geiss} for the case $n=2$. 
We consider it an idea of great significance in the theory of rigid Cohen-Macaulay modules, and it will be interesting if we can overcome the combinatorial difficulty which appears in the procedure of mutation. 


\vskip1em
Now we briefly describe the computational results obtained in this paper. 

Let $k$  be an algebraically closed field of characteristic zero, and 
let  $G$  be a finite subgroup of  $\GL (d, k)$  for  an integer  $d$. 
The group $G$ acts linearly on the formal power series ring  $S = k[[x_1,x_2, \ldots , x_d]]$  in the natural way.
We denote the invariant subring  by  $R$, i.e.  $R = S ^G$. 
It is known that  $R$  is  a Cohen-Macaulay ring. 
We are interested in maximal Cohen-Macaulay modules over  $R$  which are rigid. 
\begin{defn}
An $R$-module  $M$ is called {\it rigid} if  $\Ext^1_R (M, M) =0$. 
\end{defn}

By definition, a rigid module has no nontrivial infinitesimal deformations.  
Our aim is to classify all the rigid Cohen-Macaulay modules over $R$  up to isomorphism. 

In this paper, we accomplish a complete classification for cases involving the following two finite subgroups.

\begin{itemize}
\item[(1)]
Let  $G_1$  be a cyclic subgroup of $\GL (3, k)$  that is generated by  
$$
\sigma = 
\begin{pmatrix}
\omega && \\
&\omega & \\
&&\omega  \\
\end{pmatrix},
$$
where  $\omega$ is a primitive cubic root of unity. 
In this case,  the local ring $R = S^{G_1}$  is (the completion of) the Veronese subring of dimension three and degree three:
$$
R = k [[ \{ \text{monomials of degree three in }  x_1, x_2, x_3 \} ]]   
$$

\item[(2)]
Let  $G_2$  be a cyclic subgroup of  $\GL (4, k)$  that is generated by 
$$
\tau = 
\begin{pmatrix}
-1 &&& \\
&-1 && \\
&& -1 &\\
&&& -1 \\
\end{pmatrix}.
$$
In this case,  the local ring  $R = S ^{G_2}$  is the Veronese subring of dimension 4 and degree 2. 
$$
R = k [[ \{ \text{monomials of degree two in }  x_1, x_2, x_3, x_4 \} ]]   
$$
\end{itemize}

In both cases, it is easy to verify that  $R$ is a Gorenstein local ring with only an isolated singularity.

Our computational results are the following.

\begin{thm}\label{main1}
Let  $G = G_1$ as in case  $(1)$ above. 
Write  $S$  as a sum of modules of semi-invariants 
$S = S_0 \oplus S_1 \oplus  S_2$. Set  $M_{2i}=\Omega_R^iS_1$  and  $M_{2i+1}=\Omega_R^iS_2$ for any $i\in\Z$. 
Then a maximal Cohen-Macaulay $R$-module is rigid if and only if it is isomorphic to $R^a \oplus M_i^b \oplus M_{i+1}^c$ for some $i\in\Z$ and $a, b, c\in\Z_{\ge0}$.
\end{thm}

\begin{thm}\label{main2}
Let  $G = G_2$ as in case  $(2)$ above, and  write  $S$  as a sum of 
modules of semi-invariants  $S = S_0 \oplus S_1$.  
Then a maximal Cohen-Macaulay $R$-module is rigid if and only if it is isomorphic to $R^a\oplus(\Omega^i_RS_1)^b$ for some $i\in\Z$ and $a,b\in\Z_{\ge0}$.
\end{thm}

The rest of the paper is organized as follows. We give proofs of the theorems in Section 9. 
The other sections are devoted to theoretical preparations for the proofs and to providing the background.
In Section 2, we introduce the idea of torsion theory for a triangulated category, which simplifies the proofs in later sections. 
We also introduce the idea of mutation pairs in triangulated categories and study their basic properties.
In Section 3, the ideas of $n$-cluster tilting subcategories and AR $(n+2)$-angles are introduced and we discuss some of their elementary properties. In Section 4, we show that certain subfactor categories of triangulated categories again form triangulated categories. We also study the relationships between these triangulated categories. In particular, we find a one-one correspondence between their $n$-cluster tilting subcategories.
The idea of mutation of $n$-cluster tilting subcategories is introduced in Section 5, where we discuss some properties of mutation. 
As was mentioned above, we believe this idea is a fundamental tool for the classification of rigid objects in a triangulated category. 
We also study how many $n$-cluster tilting subcategories contain an `almost $n$-cluster tilting' subcategory.
In Section 6, we show that, in several cases, a subfactor category of a triangulated category will be equivalent to a full subcategory of an abelian category. 
 This type of argument together with Kac's theorem discussed in Section 7 plays an essential role in the proof given in Section 9. 
In Section 8, we discuss the properties of Gorenstein quotient singularities that will be necessary for our proof.
In Section 9, the properties derived for abstract triangulated categories in the preceding sections are applied to the stable category of Cohen-Macaulay modules over a Gorenstein quotient singularity, leading to the proofs of the theorems given above.  
In Section 10, we give some examples of triangulated categories based on non-commutative algebras such that we can classify all rigid objects and $n$-cluster tilting subcategories. 

We remark that Keller and Reiten gave another proof of \ref{main1} in \cite{KR2} with a different method.

\vskip1em
\begin{center}
{\bf Acknowledgement}
\end{center}

The authors would like to thank organizers and participants of the conference ``Representation Theory of Finite-Dimensional Algebras'' (Oberwolfach, February 2005), where parts of the results in this paper were presented \cite{Iyama3, Yoshino2}.
Part of this work was done while the first author visited Paderborn in February 2006. He would like to thank Henning Krause and people in Paderborn for hospitality and stimulating discussions.


\vskip1em
\begin{center}
{\bf Conventions}
\end{center}

All additive categories considered in this paper are assumed to be Krull-Schmidt,
i.e. any object is isomorphic to a finite direct sum of objects whose endomorphism rings are local.
Let $\cc$ be an additive category.
We denote by $\ind\cc$ the set of isoclasses of indecomposable objects in $\cc$, and by $\cc (X,Y)$ or $(X,Y)$ the set of morphisms from $X$ to $Y$ in $\cc$.
We denote the composition of $f \in \cc(X,Y)$ and $g \in \cc(Y,Z)$  by $fg \in \cc(X,Z)$.
An {\it ideal} $I$ of $\cc$ is an additive subgroup $I(X,Y)$ of $\cc(X,Y)$ such that $fgh\in I(W,Z)$ for any $f\in\cc(W,X)$, $g\in I(X,Y)$ and $h\in\cc(Y,Z)$.
We denote by $J_{\cc}$ the {\it Jacobson radical} of $\cc$. Namely, $J_{\cc}$ is an ideal of  $\cc$  such that  $J_{\cc}(X,X)$ coincides with the Jacobson radical of the endomorphism ring $\End_{\cc}(X)$ for any $X\in\cc$.
For an ideal $I$ of $\cc$, we denote by $\cc /I$  the category whose objects are objects of  $\cc$ and whose morphisms are elements of
$$
\cc (M, N) /I (M, N)  \ \ \text{for} \ \ M, N \in \cc/I.
$$

When we say that $\dd$ is a subcategory of $\cc$, we always mean that $\dd$ is a full subcategory which  is  closed under isomorphisms, direct sums and direct summands. For an object $X\in\cc$, we denote by $\add X$ the smallest subcategory of $\cc$ containing $X$. We denote by $[\dd]$ the ideal of $\cc$ consisting of morphisms which factor through objects in  $\dd$. Thus we have the category $\cc /[\dd]$. We simply denote  $\cc/[\add X]$  by $\cc/[X]$.

Let $f\in\cc(X,Y)$ be a morphism.
We call $f$ {\it right minimal} if it does not have a direct summand of the form $T\to 0$ ($T\in\cc$, $T\neq 0$) as a complex. Obviously, any morphism is a direct sum of a right minimal morphism and a complex of the form $T\to 0$ ($T\in\cc$). Dually, a {\it left minimal morphism} is defined.
We call $f$ a {\it sink map} of $Y\in\cc$ if $f$ is right minimal, $f\in J_{\cc}$ and
$$
\begin{CD}
\cc(-,X) @>\cdot f>> J_{\cc}(-,Y) @>>>0
\end{CD}
$$
is exact as functors on $\cc$. Dually, a {\it source map} is defined.
For a subcategory $\dd$ of $\cc$, we call $f$ a {\it right $\dd$-approximation} of $Y\in\cc$ \cite{AS} if $X\in\dd$ and
$$
\begin{CD}
\cc(-,X) @>\cdot f>> \cc(-,Y) @>>>0
\end{CD}
$$ is exact as functors on $\dd$. We call a right $\dd$-approximation {\it minimal} if it is right minimal.
We call $\dd$ a {\it contravariantly finite subcategory} of $\cc$ if any $Y\in\cc$ has a right $\dd$-approximation.
Dually, a {\it (minimal) left $\dd$-approximation} and a {\it covariantly finite subcategory} are defined. A contravariantly and covariantly finite subcategory is called {\it functorially finite}.

When we say that $\tt$  is a triangulated category, we always assume that $\tt$ is $k$-linear for a fixed field $k$ and $\dim_k\tt(X,Y)<\infty$ for any $X,Y \in \tt$.
Notice that we do not assume that a subcategory of $\tt$ is closed under the shift functor $[1]$ in $\tt$.

All modules over rings are left modules, and all morphisms act on modules from the right.

\section{Torsion theories and mutation pairs}

Throughout this section, let $\tt$ be a triangulated category. We introduce basic notions which will be used in this paper. Let $\xx$ and $\yy$ be subcategories of $\tt$. We put
$$
\xx^\perp := \{ T \in \tt \ |\ (\xx,T) =0 \}\ \ \ \mbox{and}\ \ \ {}^\perp\xx := \{ T \in \tt \ |\ (T,\xx)=0\}.
$$
We denote by $\xx\vee\yy$ the smallest subcategory of $\tt$ containing $\xx$ and $\yy$.
We denote by $\xx * \yy$  the collection of objects in  $\tt$ consisting of all such  $T\in\tt$  with  triangles
\begin{eqnarray}\label{triangle}
\begin{CD}
X @>a>> T @>b>> Y @>c>> X[1]  \qquad (X \in \xx, \ \ Y \in \yy).
\end{CD}
\end{eqnarray}
By the octahedral axiom, we have $(\xx*\yy)*\zz=\xx*(\yy*\zz)$. We call $\xx$ {\it extension closed} if $\xx*\xx=\xx$. Although $\xx*\yy$ is not necessarily closed under direct summands in general, we have the following sufficient conditions.

\begin{prop}\label{direct summands}
\begin{itemize}
\item[(1)] If $(\xx,\yy)=0$, then $\xx*\yy$ is closed under direct summands.

\item[(2)] If $J_{\tt}(\yy,\xx[1])=0$, then $\xx*\yy=\xx\vee\yy$. More generally, if $J_{\tt}(\yy_i,\xx_j[1])=0$ for $i\neq j$, then $(\xx_1\vee\xx_2)*(\yy_1\vee\yy_2)=(\xx_1*\yy_1)\vee(\xx_2*\yy_2)$.
\end{itemize}
\end{prop}

\begin{pf}
(1) Take a triangle (\ref{triangle}) and a decomposition $T=T_1\oplus T_2$. For a projection $p_i\in(T,T_i)$, we have a right $\xx$-approximation $ap_i\in(X,T_i)$. Thus we can take a minimal right $\xx$-approximation $a_i\in(X_i,T_i)$ and a triangle
$$
\begin{CD}
X_i @>a_i>> T_i @>>> U_i @>>> X_i[1].
\end{CD}
$$
Since $a$ is a right $\xx$-approximation, it is isomorphic to a direct sum of $a_1$, $a_2$ and a complex $(X_3\to0)$ with $X_3\in\xx$. Thus we have $Y\simeq U_1\oplus U_2\oplus X_3[1]$. This imples $U_i\in\yy$ and $T_i\in\xx*\yy$.

(2) Since the former assertion follows from the latter by putting $\xx_2=\yy_1=0$, we only show the latter half.
For  $T \in (\xx_1\vee\xx_2)*(\yy_1\vee\yy_2)$, take a triangle
$$
\begin{CD}
X_1\oplus X_2 @>>> T @>>> Y_1\oplus Y_2 @>a>> X_1[1]\oplus X_2[1]
\end{CD}
$$
with $X_i\in\xx$ and $Y_j\in\yy$.
By our assumption, $a$ is a direct sum of $a_1\in J_{\tt}(Y_1,X_1[1])$, $a_2\in J_{\tt}(Y_2,X_2[1])$ and $1_{T'} \in(T',T')$ with $T' \in\tt$.
Take a triangle
$$
\begin{CD}
X_i @>>> T_i @>>> Y_i @>a_i>> X_i[1].
\end{CD}
$$
Then we have $T_i\in\xx_i*\yy_i$ and $T\simeq T_1\oplus T_2$.
\qed
\end{pf}

\begin{defn}
We call a pair $(\xx,\yy)$ of subcategories of $\tt$ a {\it torsion theory} if
$$
(\xx,\yy)=0\ \ \ \mbox{and}\ \ \ \tt=\xx*\yy.
$$
In this case, it is easy to see that $\xx={}^\perp\yy$ and $\yy=\xx^\perp$ hold. Moreover, we see that $\xx$ (resp. $\yy$) is a contravariantly (resp. covariantly) finite and extension closed subcategory of $\tt$.
We do not assume that $\xx$ (resp. $\yy$) is closed under $[1]$ (resp. $[-1]$), so $(\xx,\yy[1])$ does not necessarily form a t-structure.
\end{defn}

We have the following Wakamatsu-type Lemma (1) and Auslander-Reiten-type correspondence (2). They are triangulated analogies of \cite{AR2.5}, and closely related to work of Keller and Vossieck \cite{KV2} and Beligiannis and Reiten \cite{BR}.

\begin{prop}\label{wakamatsu}
\begin{itemize}
\item[(1)] Let $\xx$ be a contravariantly finite and extension closed subcategory of $\tt$. Then $(\xx,\xx^\perp)$ forms a torsion theory.

\item[(2)] $\xx\mapsto\yy:=\xx^\perp$ gives a one-one correspondence between contravariantly finite and extension closed subcategories $\xx$ of $\tt$ and covariantly finite and extension closed subcategories $\yy$ of $\tt$. The inverse is given by $\yy\mapsto\xx:={}^\perp\yy$.
\end{itemize}
\end{prop}

\begin{pf}
(1) For any $T\in\tt$, take a triangle
\begin{eqnarray*}
\begin{CD}
X @>a>> T @>b>> Y @>c>> X[1]
\end{CD}
\end{eqnarray*}
with a minimal right $\xx$-approximation $a$ of $T$. Then $c\in J_{\tt}$ holds. We only have to show $Y\in\xx^\perp$. Take any morphism $f\in(X',Y)$ with $X'\in\xx$. By the octahedral axiom, we have a commutative diagram
$$
\begin{CD}
@. Z @=  Z @. \\
@. @VVV @VVV  @. \\
X @>>> X'' @>>> X' @>>> X[1] \\
@| @VdVV @VfVV @| \\
X @>a>> T @>b>> Y @>c>> X[1] \\
@. @VeVV @VgVV @. \\
@. Z @=  Z
\end{CD}
$$
of triangles. Since $X''\in\xx*\xx=\xx$, we have that $d$ factors through $a$. Thus we have $db=0$. Hence there exists $b'\in(Z,Y)$ such that $b=eb'$. Since $b(1_Y-gb')=0$, $1_Y-gb'$ factors through $c\in J_{\tt}$. Thus $g$ is a split monomorphism, and we have $f=0$.

(2) The assertion follows from (1) because $\xx^\perp$ is a covariantly finite and extension closed subcategory of $\tt$ satisfying ${}^\perp(\xx^\perp)=\xx$.
\qed
\end{pf}

\begin{prop}\label{torsion}
Let $\cc_i$ be a contravariantly finite and extension closed subcategory of $\tt$ such that $(\cc_i,\cc_j[1])=0$ for any $i<j$. Put
$$
\xx_n:=\add(\cc_1*\cc_2*\cdots*\cc_n),\ \ \ \ \ \yy_n:=\bigcap_{i=1}^n\cc_i^\perp.$$
Then $(\xx_n,\yy_n)$ forms a torsion theory.
\end{prop}

\begin{pf}
The case $n=1$ is proved in \ref{wakamatsu}.
We show $\tt=\xx_n*\yy_n$ by induction on  $n$.
Assume that the assertion is true for $n=i-1$. For any $T\in\tt$, take a triangle
$$
\begin{CD}
X_{i-1} @>>> T @>>> Y_{i-1} @>>> X_{i-1}[1]
\end{CD}
$$
with $X_{i-1}\in\xx_{i-1}$ and $Y_{i-1}\in\yy_{i-1}$. Take a triangle
$$
\begin{CD}
C_i @>a_i>> Y_{i-1} @>>> Y_{i} @>>> C_i[1]
\end{CD}
$$
with a minimal right $\cc_i$-approximation $a_i$ of $Y_{i-1}$. By \ref{wakamatsu}, $Y_i\in\cc_i^\perp$. Since both of $Y_{i-1}$ and $C_i[1]$ belong to $\yy_{i-1}$, we have $Y_i\in\cc_i^\perp\cap\yy_{i-1}=\yy_i$. By the octahedral axiom, we have a commutative diagram
$$
\begin{CD}
@. Y_i[-1] @=  Y_i[-1] @. \\
@. @VVV @VVV  @. \\
X_{i-1} @>>> X_i @>>> C_i @>>> X_{i-1}[1] \\
@| @VVV @Va_iVV @| \\
X_{i-1} @>>> T @>>> Y_{i-1} @>>> X_{i-1}[1]. \\
@. @VVV @VVV @. \\
@. Y_i @=  Y_i
\end{CD}
$$
of triangles. Since $X_i\in\xx_{i-1}*\cc_i=\xx_i$, we have $T\in\xx_i*\yy_i$ by the left vertical triangle.
\qed
\end{pf}





Now we consider a correspondence of subcategories of $\tt$.
This includes the mutation of cluster tilting objects and maximal 1-orthogonal objects given by Buan, Marsh, Reineke, Reiten and Todorov \cite{Buan} and Geiss, Leclerc and Schr\"oer \cite{Geiss} respectively.
In section 5, we shall apply this to introduce the mutation of $n$-cluster tilting subcategories.

\begin{defn}\label{mutationpair}
Fix a subcategory $\dd$ of $\tt$ satisfying $(\dd,\dd[1])=0$. For
a subcategory $\xx$ of $\tt$, put
$$
\mu^{-1}(\xx;\dd):=(\dd*\xx[1])\cap{}^\perp\dd[1].
$$
Then $\mu^{-1}(\xx;\dd)$ consists of all $T\in\tt$ such that there
exists a triangle
$$
\begin{CD}
X @>a>> D @>>> T @>>> X[1]
\end{CD}
$$
with $X\in\xx$ and a left $\dd$-approximation $a$. Dually, for a
subcategory $\yy$ of $\tt$, put
$$
\mu(\yy;\dd):=(\yy[-1]*\dd)\cap\dd[-1]^\perp.
$$
Then $\mu(\yy;\dd)$ consists of all $T\in\tt$ such that there
exists a triangle
$$
\begin{CD}
T @>>> D @>b>> Y @>>> T[1]
\end{CD}
$$
with $Y\in\yy$ and a right $\dd$-approximation $b$.

We call a pair $(\xx,\yy)$ of subcategories of $\tt$ a {\it
$\dd$-mutation pair} if
$$
\dd\subset\yy\subset\mu^{-1}(\xx;\dd)\ \ \ \mbox{and}\ \ \ \dd\subset\xx\subset\mu(\yy;\dd).
$$
In this case $\dd$ is a covariantly finite subcategory of $\xx$ and a contravariantly finite subcategory of $\yy$.
For a $\dd$-mutation pair $(\xx,\yy)$, we construct a functor $\ggg:\xx/[\dd]\to\yy/[\dd]$ as follows:
For any  $X \in \xx$, fix a triangle
$$
\begin{CD}
X @>{\alpha_X}>> D_X @>{\beta_X}>> \ggg X @>{\gamma_X}>> X[1]
\end{CD}
$$
with  $D_X \in \dd$  and  $\ggg X \in \yy$, and define  $\ggg X$  by this.
Then $\alpha_X$ is a left $\dd$-approximation and $\beta_X$ is a right $\dd$-approximation.
For any morphism $f\in(X,X')$,
there exists $g$ and $h$ which make the following diagram commutative.
$$
\begin{CD}
X @>{\alpha_X}>> D_X @>{\beta_X}>> \ggg X @>{\gamma_X}>> X[1] \\
@V{f}VV  @V{g}VV  @V{h}VV @V{f[1]}VV   \\
X' @>{\alpha_{X'}}>> D_{X'} @>{\beta_{X'}}>> \ggg X' @>{\gamma_{X'}}>> X'[1]
\end{CD}
$$
Now put $\ggg\overline{f}:=\overline{h}$.
\end{defn}

\begin{prop}\label{mutationequivalence}
Under the circumstances as above, the following assertions hold.
\begin{itemize}
\item[(1)] $\ggg:\xx/[\dd]\to\yy/[\dd]$ gives a well-defined equivalence.

\item[(2)] $\yy=\mu^{-1}(\xx;\dd)$ and $\xx=\mu(\yy;\dd)$ hold.
\end{itemize}
\end{prop}

\begin{pf}
(1) Assume that the diagram
$$
\begin{CD}
X @>{\alpha_X}>> D_X @>{\beta_X}>> \ggg X @>{\gamma_X}>> X[1] \\
@V{f}VV  @V{g'}VV  @V{h'}VV @V{f[1]}VV   \\
X' @>{\alpha_{X'}}>> D_{X'} @>{\beta_{X'}}>> \ggg X' @>{\gamma_{X'}}>> X'[1]
\end{CD}
$$
is also commutative.
Then, since $(h'-h)\gamma_{X'}=0$ holds,
$h'-h$ factors through $\beta_{X'}$.
Thus $\overline{h'}=\overline{h}$ holds.
This shows that $\ggg$ is a well-defined functor.

For any $Y\in\yy$, we fix a triangle
$$
\begin{CD}
Y[-1] @>>> \hhh Y @>>> D^Y @>>> Y
\end{CD}
$$
with $D^Y\in\dd$ and $\hhh Y\in\xx$.
We can construct a functor $\hhh:\yy/[\dd]\to\xx/[\dd]$ in a dual manner.
We can easily show that $\hhh$ gives a quasi-inverse of $\ggg$.

(2) For any $Y\in\mu^{-1}(\xx;\dd)\subset\dd*\xx[1]$, take a triangle
$$
\begin{CD}
X @>a>> D @>>> Y @>>> X[1]
\end{CD}
$$
with $X\in\xx$ and $D\in\dd$. Since $(Y[-1],\dd)=0$ holds, $a$ is a left $\dd$-approximation. Thus $Y$ and $\ggg X\in\yy$ are isomorphic up to a direct summand in $\dd$. Since $\yy$ contains $\dd$, we have $Y\in\yy$.
\qed
\end{pf}

\begin{prop}
Let $\dd$ be a functorially finite subcategory of $\tt$ satisfying $(\dd,\dd[1])=0$.
Then $\xx\mapsto\yy:=\mu^{-1}(\xx;\dd)$ gives a one-one correspondence between subcategories $\xx$ of $\tt$ satisfying $\dd\subset\xx\subset\dd[-1]^\perp$ and subcategories $\yy$ of $\tt$ satisfying $\dd\subset\yy\subset{}^\perp\dd[1]$. The inverse is given by $\yy\mapsto\xx:=\mu(\yy;\dd)$.
\end{prop}

\begin{pf}
For a given $\xx$, put $\yy:=\mu^{-1}(\xx;\dd)$. Then $\yy$ is a
subcategory of $\tt$ since $\dd*\xx[1]$ is closed under direct
summands by \ref{direct summands}(1). For any $X\in\xx$, take a
triangle
$$
\begin{CD}
X @>a>> D @>b>> Y @>>> X[1]
\end{CD}
$$
with a left $\dd$-approximation $a$. Then $Y\in\yy$ and $b$ is a
right $\dd$-approximation. Thus $\xx\subset\mu(\yy;\dd)$ holds, so
$(\xx,\yy)$ is a $\dd$-mutation pair. By
\ref{mutationequivalence}, we have $\xx=\mu(\yy;\dd)$.
\qed

\end{pf}

In the rest of this section, we make several preliminaries on a Serre functor on a $k$-linear triangulated category $\tt$.

\begin{defn}
Following Bondal and Kapranov \cite{Bondal1}, we call a $k$-linear autofunctor $\sss :\tt\to\tt$ a {\it Serre functor} of $\tt$ if there exists a functorial isomorphism
$$
(X,Y) \simeq D(Y,\sss X)
$$
for any $X,Y\in\tt$, where  $D$ denotes the $k$-dual. Obviously $\xx^\perp={}^\perp(\sss \xx)$ holds for a subcategory $\xx$ of $\tt$. We say that $\tt$ is  {\it $n$-Calabi-Yau} ($n\in\Z$) if $\sss=[n]$.
\end{defn}


The existence of a Serre functor is related with the classical concept of coherent functors and dualizing $k$-varieties, which were introduced by Auslander and Reiten \cite{coherent, AR} and has given a theoretical background of subsequent Auslander-Reiten theory.
Recently Krause (e.g. \cite{Krause1, Krause2}) effectively applied them to study Brown representability and the telescope conjecture.

\begin{defn}\label{coherent}
Let $\cc$ be a $k$-linear category.
A $\cc$-module is a contravariant $k$-linear functor $F:\cc\to\Mod k$.
Then $\cc$-modules form an abelian category $\Mod\cc$.
By Yoneda's lemma, representable functors are projective objects in $\Mod \cc$. The $k$-dual $D$ induces a functor  $D: \Mod \cc \leftrightarrow \Mod \cc^{\rm op}$.

We call $F \in \Mod \cc$ {\it coherent} \cite{coherent} if there exists an exact sequence
\begin{eqnarray}\label{projective resolution}
\begin{CD}
(-,Y) @>\cdot a>> (-,X) @>>> F @>>> 0
\end{CD}
\end{eqnarray}
of $\cc$-modules with $X,Y\in\tt$.
We denote by $\mod \cc$ the full subcategory of  $\Mod \cc$ consisting of coherent $\cc$-modules. It is easily checked that $\mod\cc$ is closed under cokernels and extensions in $\Mod\cc$. Moreover, $\mod\cc$ is closed under kernels in $\Mod\cc$ if and only if $\cc$ has pseudokernels. In this case, $\mod\cc$ forms an abelian category (see \cite{coherent}). For example, if $\cc$ is a contravariantly finite subcategory of a triangulated category, then $\cc$ has pseudokernels and $\mod\cc$ forms an abelian category.

We call $\cc$ a {\it dualizing $k$-variety} \cite{AR} if $D$ gives a duality $\mod \cc \leftrightarrow \mod \cc^{\rm op}$. Then $\cc$ has pseudokernels and $\mod\cc$ is an abelian subcategory of $\Mod\cc$. In this case, for any $F\in\mod\cc$, take a minimal projective resolution (\ref{projective resolution}) and define the {\it Auslander-Reiten translation} $\tau F\in\mod\cc$ by an exact sequence
$$
\begin{CD}
0 @>>> \tau F @>>> D(Y,-) @>D(a\cdot)>> D(X,-).
\end{CD}
$$
We collect some properties of dualizing $k$-varieties.
\end{defn}

\begin{prop}\label{dualizing}
\begin{itemize}
\item[(1)] Any functorially finite subcategory of a dualizing $k$-variety is also a dualizing $k$-variety.

\item[(2)] Any object in a dualizing $k$-variety has a sink map and a source map.
\end{itemize}
\end{prop}

\begin{pf}
Although each assertion is well-known (see \cite{AR, AS}), we shall give a proof of (2) for the convenience of the reader.


(2) In general, it follows immediately from the definition that $X\in\cc$ has a sink map if and only if the simple $\cc$-module $S_X:=(\cc/J_{\cc})(-,X)$ is coherent. We have epimorphisms $f:(-,X)\to S_X$ and $g:(X,-)\to DS_X$. Then $S_X$ is an image of the morphism $f(Dg):(-,X)\to D(X,-)$ in $\mod\cc$. Since $\mod\cc$ is closed under images in $\Mod\cc$, we have $S_X\in\mod\cc$.
\qed
\end{pf}

We now have a relationship between Serre functors and dualizing $k$-varieties (cf. \cite[2.2]{AR3}).


\begin{prop}\label{fun}
A triangulated category  $\tt$ has a Serre functor if and only if it is a dualizing $k$-variety.
\end{prop}

\begin{pf}
To prove the \lq only if' part, suppose the triangulated category  $\tt$  has a Serre functor $\sss$.
Since $D(-,T)\simeq(\sss^{-1}T,-)$ holds for any $T\in\tt$, we have $DF\in\mod\tt^{\rm op}$ for any representable functor $F$.
Since $\tt^{\rm op}$ is a triangulated category, $\mod\tt^{\rm op}$ is closed under kernels. Thus we have $DF\in\mod\tt^{\rm op}$ for any $F\in\mod\tt$.


We show the \lq if'  part. Fix  $T \in \tt$. Since $D (T,-) \in \mod \tt$,  there is an exact sequence $(-,Y) \stackrel{\cdot a}{\to} (-,X) \to D(T,-) \to 0$. Take a triangle involving $a$ as $Z \stackrel{b}{\to} Y \stackrel{a}{\to} X \stackrel{c}{\to} Z[1]$ in $\tt$. Then we have an exact sequence $0 \to D(T,-) \to (-,Z[1]) \stackrel{\cdot b[1]}{\to} (-,Y[1])$, which splits since $D(T,-)$ is an injective object in $\mod\tt$. Thus $D(T,-)$ is representable, so there exists an object $\sss T \in \tt$ such that $D(T,-) \simeq (-,\sss T)$. One easily verifies that $\sss$ gives a Serre functor of $\tt$.
\qed
\end{pf}


\section{$n$-cluster tilting subcategories}

Let $\tt$ be a triangulated category. For an integer  $n>0$, we call a subcategory $\cc$ of $\tt$ {\it $n$-rigid} if the equalities
$$
(\cc,\cc[i])=0 \qquad (0 < i <n)
$$
hold. We often call a 2-rigid subcategory {\it rigid}. More strongly, a functorially finite subcategory $\cc$ of $\tt$ is said to be {\it $n$-cluster tilting} \cite{KR} if it satisfies
$$
\cc=\bigcap_{i=1}^{n-1} \cc[-i]^\perp = \bigcap_{i=1}^{n-1}{}^\perp\cc[i].
$$
This is a triangulated analog of `maximal $(n-1)$-orthogonal subcategories' studied in \cite{Iyama1, Iyama2, Iyama3}.
Of course, $\tt$ is a unique $1$-cluster tilting subcategory of $\tt$. It is obvious that if $\dd$ is an $n$-rigid subcategory containing an $n$-cluster tilting subcategory $\cc$, then it holds that $\cc = \dd$.

\begin{thm}\label{n-cluster torsion}
Let $\cc$ be an $n$-cluster tilting subcategory of $\tt$.
\begin{itemize}
\item[(1)] $\tt=\cc*\cc[1]*\cdots*\cc[n-1]$ holds.

\item[(2)] $(\cc*\cc[1]*\cdots*\cc[\ell-1],\ \cc[\ell]*\cc[\ell+1]*\cdots*\cc[n-1])$ forms a torsion theory for any $\ell$ ($0<\ell<n$).
\end{itemize}
\end{thm}

\begin{pf}
Since $\cc$ is $n$-cluster tilting, we have $\displaystyle\bigcap_{i=0}^{n-2}\cc[i]^\perp=\cc[n-1]$. By \ref{direct summands} and \ref{torsion}, we have a torsion theory $(\cc*\cc[1]*\cdots*\cc[n-2],\ \cc[n-1])$. Thus (1) follows. Now (2) follows from (1) and \ref{direct summands} because we have $(\cc*\cdots*\cc[\ell-1],\ \cc[\ell]*\cdots*\cc[n-1])=0$.
\qed
\end{pf}

\begin{notation}\label{1.6}
Let $n\ge 1$ be an integer.
Suppose we are given triangles in $\tt$.
$$
\begin{CD}
X_{i+1} @>{b_{i+1}}>> C_i @>{a_i}>> X_i @>>> X_{i+1}[1] \qquad (0\le i<n)
\end{CD}
$$
Then we call a complex
$$
\begin{CD}
X_{n} \stackrel{b_{n}}{\longrightarrow} C_{n-1} \stackrel{a_{n-1}b_{n-1}}{\longrightarrow} C_{n-2} \stackrel{a_{n-2}b_{n-2}}{\longrightarrow} @.\cdots\cdots\cdots @. \stackrel{a_{2}b_{2}}{\longrightarrow} C_{1} \stackrel{a_{1}b_{1}}{\longrightarrow} C_0 \stackrel{a_0}{\longrightarrow} X_0
\end{CD}
$$
an {\it $(n+2)$-angle} in $\tt$.
We sometimes denote it by
\begin{eqnarray}\label{(n+2)-angle}
\begin{CD}
X_{n} \stackrel{b_{n}}{\longrightarrow} C_{n-1}\ang{X_{n-1}}{}{a_{n-1}}{b_{n-1}}C_{n-2}\ang{X_{n-2}}{}{a_{n-2}}{b_{n-2}} @. \cdots\cdots\cdots @. \ang{X_2}{}{\ \ a_2}{b_2}C_1\ang{X_1}{}{\ \ a_1}{b_1}C_0 \stackrel{a_0}{\longrightarrow} X_0.
\end{CD}
\end{eqnarray}
\end{notation}

The following corollary is a triangulated version of \cite[2.2.3]{Iyama1}.

\begin{cor}\label{1.7}
Let $\cc$ be an $n$-cluster tilting subcategory of $\tt$ and $X_0\in\tt$. Then there exists an $(n+2)$-angle
$$
0\stackrel{}{\longrightarrow}C_{n-1}\stackrel{f_{n-1}}{\longrightarrow}C_{n-2}\stackrel{f_{n-2}}{\longrightarrow}\cdots\cdots\stackrel{f_2}{\longrightarrow}C_1\stackrel{f_1}{\longrightarrow}C_0\stackrel{}{\longrightarrow}X_0
$$
with $C_i\in\cc \ (0\le i<n)$ and $f_i\in J_{\cc} \ (0<i<n)$.
\end{cor}

\begin{pf}
Since $X_0\in\tt=\cc*\cdots*\cc[n-1]$ by \ref{n-cluster torsion}, we can inductively construct a triangle $X_{i+1}\stackrel{b_{i+1}}{\to} C_i\stackrel{}{\to}X_i\to X_{i+1}[1]$ with $C_i\in\cc$ and $X_i\in\cc*\cdots*\cc[n-1-i]$. We can assume $b_{i+1}\in J_{\cc}$ by \ref{direct summands}. By gluing them, we have the desired $(n+2)$-angle.
\qed
\end{pf}

In the rest of this section, we assume that $\tt$ has a Serre functor $\sss$. We put
$$
\sss_n = \sss \circ[-n] : \tt \to \tt.
$$
By definition, $\tt$ is $n$-Calabi-Yau if and only if $\sss_n$ is an identity functor.
We call a subcategory $\cc$ of $\tt$ an {\it $\sss_n$-subcategory} of $\tt$  if it satisfies
$$\cc=\sss_n\cc=\sss_n^{-1}\cc.$$

\begin{prop}\label{1.3}
Any $n$-cluster tilting subcategory $\cc$ of $\tt$ is an $\sss_n$-subcategory.
In particular, $\cc$ has an autofunctor $\sss_n:\cc\to\cc$.
\end{prop}

\begin{pf}
$\displaystyle
\cc = \bigcap_{i=1}^{n-1}\cc[-i]^\perp
= \bigcap_{i=1}^{n-1}{}^\perp \sss\cc[-i]
= \sss_n\bigcap_{i=1}^{n-1}{}^\perp\cc[n-i]=\sss_n\cc.\mbox{\qed}
$
\end{pf}

The autofunctor $\sss_n$ of $\cc$ is a triangulated analog of the $n$-Auslander-Reiten translation functor in \cite{Iyama1}.

It is easily checked that any $\sss_n$-subcategory $\cc$ satisfies $\displaystyle\displaystyle\bigcap_{i=1}^{n-1}\cc[-i]^\perp=\bigcap_{i=1}^{n-1}{}^\perp\cc[i]$. Thus we obtain the following proposition.

\begin{prop}\label{1.4}
For a subcategory $\cc$ of $\tt$, the following conditions are equivalent.
\begin{itemize}
\item[(1)]
 $\cc$ is an $n$-cluster tilting subcategory of $\tt$.

\item[(2)]
 $\cc$ is an $\sss_n$-subcategory of $\tt$ satisfying
$\displaystyle\cc=\bigcap_{i=1}^{n-1}\cc[-i]^\perp$.

\item[(3)]
 $\cc$ is an $\sss_n$-subcategory of $\tt$ satisfying
$\displaystyle\cc=\bigcap_{i=1}^{n-1}{}^\perp\cc[i]$.
\end{itemize}
\end{prop}

\begin{lemma}\label{1.8}
Let  $\cc$  be an $n$-rigid subcategory of  $\tt$, and let all terms $X_0, X_n$  and $C_i$
in an $(n+2)$-angle (\ref{(n+2)-angle}) be in $\cc$.
Define a $\cc$-module $F$ and a $\cc^{\rm op}$-module $G$ by
$$
(-,C_0)\stackrel{\cdot a_0}{\longrightarrow}(-, X_0)\longrightarrow F\longrightarrow0\ \ \ \mbox{and}\ \ \ (C_{n-1},-)\stackrel{b_n \cdot}{\longrightarrow}(X_n,-)\longrightarrow G\longrightarrow0.
$$
\begin{itemize}
\item[(1)]
For any $i \ \ (0<i<n)$,  we have the following isomorphisms of $\cc$-modules:
\begin{eqnarray*}
(-,X_i[j])&\simeq
&\begin{cases}
\ F&(j=i)\\
\ 0&(0<j<n,\ j\neq i),
\end{cases}
\\
(X_{n-i}[-j],-)&\simeq
&\begin{cases}
\ G&(j=i)\\
\ 0&(0<j<n,\ j\neq i),
\end{cases}
\end{eqnarray*}

\item[(2)]
Moreover, if $\cc$ is an $\sss_n$-subcategory, then we have
$$
G\simeq D(F\circ \sss_n^{-1})\ \ \ \mbox{and}\ \ \ F\simeq D(G\circ \sss_n).
$$
\end{itemize}
\end{lemma}

\begin{pf}
(1) First we note that there is an isomorphism
{\small
 $$
\begin{CD}
0=(-,C_i[j]) @>>> (-, X_{i}[j]) @>{\simeq}>> (-, X_{i+1}[j+1]) @>>> (-,C_{i}[j+1])=0
\end{CD}
$$
}as functors on $\cc$  if  $0< j <n-1$.
Therefore, if $0< j < i <n$, then we have isomorphisms
$$
(-,X_{i}[j])\simeq(-,X_{i+1}[j+1])\simeq\cdots\simeq(-,X_n[n+j-i])=0,
$$
since  $X_n \in \cc$  and  $0< n+j-i < n$.
On the other hand, in the case  $0< i \le j < n$,  we have the exact sequence
$$
(-,C_0[j-i])\to(-,X_0[j-i])\to(-,X_{1}[j-i+1])\to(-,C_{0}[j-i+1]),
$$
hence
{\small
$$
(-,X_{i}[j])\simeq(-,X_{i-1}[j-1])\simeq\cdots\cdots\simeq(-,X_1[j-i+1])\simeq
\begin{cases}
F&(i=j), \\
(-, X_0[j-i]) = 0 &(i<j),  \\
\end{cases}
$$
}since  $X_0 \in \cc$.
We have proved the first isomorphism, and the second isomorphism is proved in a similar manner.

(2) Consider a commutative diagram with exact rows consisting of functors on  $\cc$.
{\small
$$
\begin{CD}
0 @>>> D(F\circ \sss_n^{-1}) @>>> D(\sss_n^{-1}-,X_0)@>>> D(\sss_n^{-1}-,C_0)\\
@. @. @| @| \\
(C_0[1-n],-) @>>> (X_1[1-n],-) @>>> (X_0[-n],-) @>>> (C_0[-n],-)
\end{CD}
$$
}
If $n=1$, then we immediately obtain an isomorphism $D(F\circ \sss_n^{-1})\simeq G$.
If $n>1$, then since  $(C_0 [1-n], -) = 0$  as a functor on  $\cc$, we obtain $D(F\circ \sss_n^{-1})\simeq(X_1[1-n],-)$, which is isomorphic to $G$ by the second isomorphism in (1).
The second isomorphism follows from the first one.
\qed\end{pf}

\begin{cor}\label{(n+2)-angle approximation}
Let  $\cc$  be an $n$-rigid subcategory of  $\tt$, and let all terms $X_0, X_n$ and $C_i$
in an $(n+2)$-angle (\ref{(n+2)-angle}) be in $\cc$.
Then $a_i$ is a right $\cc$-approximation and $b_i$ is a left $\cc$-approximation for $0<i<n$.
\end{cor}

\begin{pf}
For each $i$ with $0< i < n$, there is an exact sequence of functors on $\cc$
$$
\begin{CD}
(-,  C_i)  @>\cdot a_i>> (-, X_i)  @>>> (-, X_{i+1}[1]),
\end{CD}
$$
where $(-, X_{i+1}[1]) =0$  by \ref{1.8} and by  $X_n \in \cc$.
Therefore  $a_i$  is a right $\cc$-approximation for $0< i <n$.
Similarly, $b_i$  is a right $\cc$-approximation for $0< i <n$.
\qed
\end{pf}

\begin{defn}\label{1.9}
Let $\cc$ be an $n$-cluster tilting subcategory of $\tt$.
We call an $(n+2)$-angle
\vspace{3pt}
\begin{eqnarray*}
Y\stackrel{b_n}{\longrightarrow}C_{n-1}\ang{X_{n-1}}{}{a_{n-1}}{b_{n-1}}C_{n-2}\ang{X_{n-2}}{}{a_{n-2}}{b_{n-2}}\cdots\cdots\ang{X_2}{}{\ \ a_2}{b_2}C_1\ang{X_1}{}{\ \ a_1}{b_1}C_0\stackrel{a_0}{\longrightarrow}X
\end{eqnarray*}
an {\it AR (=Auslander-Reiten) $(n+2)$-angle} if the following conditions are satisfied\footnote{Notice that the induced complex
\begin{eqnarray*}
(-,Y)\to(-,C_{n-1})\to(-,C_{n-2})\to\cdots\to(-,C_1)\to(-,C_0)\to(-,X)
\end{eqnarray*}
is not necessarily exact.
Thus the definition given in \cite[2.1(2)]{Iyama3} contains an error.}.

\begin{itemize}
\item[(0)]
$X$, $Y$  and $C_i \ (0 \leq i <n)$ belong to  $\cc$.
\item[(i)]
$a_0$ is a sink map of $X$ in $\cc$ and $b_n$ is a source map of $Y$ in $\cc$.

\item[(ii)]
$a_i$ is a minimal right $\cc$-approximation of $X_i$ for any $i$ ($0<i<n$).

\item[(iii)]
$b_i$ is a minimal left $\cc$-approximation of $X_i$ for any $i$ ($0<i<n$).
\end{itemize}
\end{defn}

An AR $(n+2)$-angle with right term $X$
 (resp. left term $Y$) depends only on $X$ (resp. $Y$) and is unique up to isomorphism
as a complex because sink (resp. source) maps and minimal right (resp. left) $\cc$-approximations are
unique up to isomorphism of complexes.

\begin{prop}\label{1.10}
Let $\cc$ be an $n$-cluster tilting subcategory of $\tt$ and let
\begin{eqnarray}\label{angle}
Y\stackrel{f_n}{\longrightarrow}C_{n-1}\stackrel{f_{n-1}}{\longrightarrow}C_{n-2}\stackrel{f_{n-2}}{\longrightarrow}\cdots\cdots\stackrel{f_2}{\longrightarrow}C_1\stackrel{f_1}{\longrightarrow}C_0\stackrel{f_0}{\longrightarrow}X
\end{eqnarray}
be an $(n+2)$-angle with all $X$, $Y$  and  $C_i$  being in $\cc$.
Then the following conditions are equivalent.

\begin{itemize}
\item[(1)]
The $(n+2)$-angle (\ref{angle}) is an AR $(n+2)$-angle.
\item[(2)]
$f_i\in J_{\cc}\ ( 0 \leq i \leq n)$, and $f_0$ is a sink map in $\cc$.
\item[(3)]
$f_i\in J_{\cc}\ ( 0 \leq i \leq n)$, and $f_n$ is a source map in $\cc$.
\end{itemize}
\end{prop}

\begin{pf}
We use the notation in \ref{1.9}, and put $X_0:=X$  and  $X_{n}:=Y$.
One can easily check that $a_{i-1}$ is right minimal if and only if $f_{i} \in J_{\cc}$ if and only if $b_{i+1}$ is left minimal. In particular, (1) implies (2) and (3).

We only prove the implication  $(2) \Rightarrow (1)$.
By \ref{(n+2)-angle approximation}, $a_i$ (resp. $b_i$) is a minimal right (resp. left) $\cc$-approximation for $0< i <n$.
Let  $F$  and  $G$  be the same as in \ref{1.8}.
Then, since  $f_0$  is a sink map,  $F$  is a semisimple functor.
It then follows from the isomorphism  $G \simeq D(F\circ \sss_n ^{-1})$ that  $G$ is also
semisimple.
This shows that  $f_n$  is a source map.
\qed
\end{pf}

Now we state the main theorem in this section, which is a triangulated analog of \cite[3.3.1]{Iyama1}.

\begin{thm}
Let $\cc$ be an $n$-cluster tilting subcategory of $\tt$.
For any $X\in\cc \ \ (\text{resp.} \ \  Y\in\cc)$, there exists an AR $(n+2)$-angle (\ref{angle})
with $Y=\sss_nX \ \ (\text{resp.} \ \ X=\sss_n^{-1}Y)$.
\end{thm}

\begin{pf}
Starting from $X\in\cc$, take a triangle
$$
\begin{CD}
X_1 @>>>  C_0 @>{f_0}>> X @>>>  X_1[1]
\end{CD}
$$
with a sink map $f_0$ in $\cc$ by \ref{dualizing} and \ref{fun}.
Then we apply  \ref{1.7}  to $X_1$ to get
 an $(n+1)$-angle
$$
C_n\stackrel{f_n}{\longrightarrow}C_{n-1}\stackrel{f_{n-1}}{\longrightarrow}C_{n-2}\stackrel{f_{n-2}}{\longrightarrow}\cdots\cdots\stackrel{f_2}{\longrightarrow}C_1\stackrel{}{\longrightarrow}X_1
$$
with $C_i\in\cc \ (0<i\le n)$ and $f_i\in J_{\cc} \ (1<i\le n)$.
By gluing them and applying \ref{1.8}, we obtain an AR $(n+2)$-angle
$$
C_n\stackrel{f_n}{\longrightarrow}C_{n-1}\stackrel{f_{n-1}}{\longrightarrow}C_{n-2}\stackrel{f_{n-2}}{\longrightarrow}\cdots\cdots\stackrel{f_2}{\longrightarrow}C_1\stackrel{f_1}{\longrightarrow}C_0\stackrel{f_0}{\longrightarrow}X.
$$
It remains to show $C_n=\sss_nX$. Let  $F$ and $G$  be the same as in \ref{1.8}. Then we have $F=(\cc/J_{\cc})(-,X)$ and
$$
G=D(F\circ\sss_n^{-1})=D(\cc/J_{\cc})(\sss_n^{-1}-,X)\simeq D(\cc/J_{\cc})(-,\sss_nX)\simeq(\cc/J_{\cc})(\sss_nX,-).
$$
Since $(C_{n-1},-)\stackrel{f_n\cdot}{\to}(C_n,-)\to G\to0$ is exact with $f_n\in J_{\cc}$, we have $C_n=\sss_nX$.
\qed
\end{pf}
\section{Subfactor triangulated categories }


Throughout this section, let $\tt$ be a triangulated category and let $\zz \supset \dd$ be subcategories of $\tt$. We assume the following two conditions concerning  $\zz$  and  $\dd$:
\begin{itemize}
\item[(Z1)]
$\zz$ is extension closed,
\item[(Z2)]
$(\zz,\zz)$ forms a $\dd$-mutation pair.
\end{itemize}
Under such a setting,  we put
$$
\uu := \zz/[\dd].
$$
The aim of this section is to prove that $\uu$ actually forms a triangulated category.

\begin{defn}\label{subfactor3.2}
Let
$$
\langle1\rangle:=\ggg:\uu\to\uu
$$
be the equivalence constructed in \ref{mutationequivalence}.
Thus for any  $X \in \zz$, we have a triangle
$$
\begin{CD}
X @>{\alpha_X}>> D_X @>{\beta_X}>> X \langle1\rangle @>{\gamma_X}>> X[1]
\end{CD}
$$
where  $D_X \in \dd$  and  $X \langle1\rangle \in \zz$.

Let
$$
\begin{CD}
X @>{a}>> Y @>{b}>> Z @>{c}>> X[1]
\end{CD}
$$
be a triangle in $\tt$ with $X,Y,Z\in\zz$.
Since $\tt(Z[-1],D_X)=0$ holds,  there is a commutative diagram of triangles :
\begin{eqnarray}\label{triangle in U}
\begin{CD}
X @>{a}>> Y @>{b}>> Z @>{c}>> X[1] \\
@|   @VVV @V{d}VV @|  \\
X @>{\alpha_X}>> D_X @>{\beta_X}>> X\langle1\rangle @>{\gamma_X}>> X[1]
\end{CD}
\end{eqnarray}
Now we consider the complex
$$
\begin{CD}
X @>{\overline{a}}>> Y @>{\overline{b}}>> Z @>{\overline{d}}>> X \langle1\rangle
\end{CD}
$$
in $\uu$.
We define {\it triangles in $\uu$} as the complexes obtained in this way.
\end{defn}

The aim of this section is to prove the following theorem.

\begin{thm}\label{subfactor3.3}
The category $\uu$ forms a triangulated category with respect to the auto-equivalence $\langle1\rangle$ and triangles defined in \ref{subfactor3.2}.
\end{thm}

Let us start with the lemma below.
We call  $a\in\tt(X,Y)$  {\it $\dd$-monic}  (resp. $\dd$-epic) if $\tt(Y,\dd)\stackrel{a\cdot}{\to}\tt(X,\dd)\to0$
(resp. $\tt(\dd,X)\stackrel{\cdot a}{\to}\tt(\dd,Y)\to0$) is exact.

\begin{lemma}\label{subfactor3.4}
\begin{itemize}
\item[(1)]
For any $X, Y \in \zz$ and any $\overline{a}\in\uu(X,Y)$, there exists $\dd$-monic $a'\in\tt(X,Y')$ such that $Y'=Y$ in $\uu$ and $\overline{a'}=\overline{a}$.

\item[(2)]
Let $X \stackrel{a}{\to} Y \stackrel{b}{\to}Z\stackrel{c}{\to}X[1]$ be a triangle in $\tt$ with $X,Y\in\zz$.
Then $Z\in\zz$ if and only if $a$ is $\dd$-monic.

\item[(3)]
Any commutative diagram
$$
\begin{CD}
X @>{a}>> Y @>{b}>> Z @>{c}>> X[1] \\
@V{f}VV  @V{g}VV @V{h}VV @V{f[1]}VV \\
X' @>{a'}>> Y' @>{b'}>> Z' @>{c'}>> X'[1]
\end{CD}
$$
of triangles in $\tt$ with $X,Y,Z,X',Y',Z'\in\zz$ induces a commutative diagram
$$
\begin{CD}
X @>{\overline{a}}>> Y @>{\overline{b}}>> Z @>{\overline{d}}>> X\langle1\rangle \\
@V{\overline{f}}VV @V{\overline{g}}VV @V{\overline{h}}VV @V{\overline{f}\langle1\rangle}VV \\
X' @>{\overline{a'}}>> Y' @>{\overline{b'}}>> Z' @>{\overline{d'}}>> X'\langle1\rangle.
\end{CD}
$$
of triangles in $\uu$.
\end{itemize}
\end{lemma}

\begin{pf}
(1)  $a' := (a\ \alpha_X) \in \tt (X,Y\oplus D_X) $  satisfies the conditions.

(2)  The `only if' part follows by taking $\tt(-,\dd)$ and using $\tt(Z[-1],\dd)=0$.
We show the `if' part.
Since $Z \in \zz*\zz[1]  \subset \zz*\zz*\dd[1]  \subset \zz*\dd[1]$  by the conditions  (Z1)  and  (Z2),
we have a triangle
$W \to Z \stackrel{d}{\to} D[1] \to W[1]$
with $W\in\zz$ and $D\in\dd$.
Since $bd=0$ by $\tt(Y,D[1])=0$, we can write $d=cd'$ with $d'\in\tt(X[1],D[1])$.
Since $a[1]$ is $\dd[1]$-monic, $d'$ factors through $a[1]$.
Thus we have $d=0$ and $Z\in\zz$.

(3) By the construction in \ref{mutationpair},
we have a commutative diagram
$$
\begin{CD}
Z @>{d}>> X \langle1\rangle @>{\gamma_X}>> X[1]\\
@. @V{f'}VV @V{f[1]}VV  \\
Z' @>{d'}>> X'\langle1\rangle @>{\gamma_{X'}}>> X'[1]
\end{CD}
$$
in $\tt$ with $\overline{f'} = \overline{f} \langle1\rangle$,
$d\gamma_X=c$ and $d'\gamma_{X'}=c'$.
Since
$$
(df'-hd')\gamma_{X'}=d\gamma_Xf[1]-hc'=cf[1]-hc'=0
$$
holds,  $df'-hd'$ factors through $\dd$.
Thus we obtain a desired commutative diagram.
\qed
\end{pf}

\vskip.5em
\noindent{\bf Proof of \ref{subfactor3.3}}
We will check the axioms of triangulated categories.

(TR0) The commutative diagram
$$
\begin{CD}
X @>{1}>> X @>>> 0 @>>> X[1] \\
@| @VVV @VVV @| \\
X @>{\alpha_X}>> D_X @>{\beta_X}>> X \langle1\rangle @>{\gamma_X}>> X[1]
\end{CD}
$$
shows that $X \stackrel{1}{\to} X \to 0 \to X \langle1\rangle$ is a triangle.

(TR1)  Fix any $\overline{a}\in\uu(X,Y)$.
By \ref{subfactor3.4}(1),  we may assume that $a\in\tt(X,Y)$ is $\dd$-monic.
By \ref{subfactor3.4}(2),  there exists a triangle $X \stackrel{a}{\to} Y \to Z \to X[1]$
with $Z\in\zz$.
Thus we have a triangle $X \stackrel{\overline{a}}{\to} Y \to Z \to X\langle1\rangle$ in $\uu$.

(TR2) Let
$X \stackrel{\overline{a}}{\longrightarrow} Y\stackrel{\overline{b}}{\longrightarrow} Z \stackrel{\overline{d}}{\longrightarrow} X \langle1\rangle$
be a triangle in $\uu$.
We may assume that it is induced by a commutative diagram (\ref{triangle in U}) of triangles in $\tt$ in \ref{subfactor3.2}.
By the octahedral axiom and $\tt(Z,D_X[1])=0$, we have a commutative diagram
$$
\begin{CD}
@. D_X @=  D_X @. \\
@. @V{(0\ 1)}VV @V{\beta_X}VV  @. \\
Y @>{b'}>> Z \oplus D_X @>{d'}>> X \langle1\rangle @>{-\gamma_X\cdot a[1]}>> Y[1] \\
@| @V{{1\choose0}}VV @V{\gamma_X}VV @| \\
Y @>{b}>> Z @>{c}>> X[1] @>{-a[1]}>> Y[1] \\
@. @V0VV @V{\alpha_X[1]}VV @. \\
@. D_X[1] @=  D_X[1]
\end{CD}
$$
of triangles in $\tt$.
Put $b'=(b_1'\ b_2')$ and  $d^{\prime}={d'_1\choose d'_2}$.
Since $(b_1'\ b_2'){1\choose0}=b$ holds,
we have $b_1'=b$ and $\overline{b'} = \overline{b}$.
Since ${d'_1\choose d'_2}\gamma_X = {1\choose 0}c$ holds,
we have $(d-d'_1)\gamma_X=0$.
Thus $d-d'_1$ factors through $\beta_X$, and we have $\overline{d'}=\overline{d}$.
Take a commutative diagram
$$
\begin{CD}
Y @>{b'}>> Z \oplus D_X @>{d'}>> X \langle1\rangle @>{-\gamma_X\cdot a[1]}>> Y[1] \\
@| @VVV @V{-a'}VV @| \\
Y @>{\alpha_Y}>> D_Y @>{\beta_Y}>> Y\langle1\rangle @>{\gamma_Y}>> Y[1]
\end{CD}
$$
of triangles in $\tt$.
Since $\gamma_X\cdot a[1]=a'\gamma_Y$ holds,
we have $\overline{a'}=\overline{a}\langle1\rangle$.
Thus we have a triangle
$$
\begin{CD}
Y @>{\overline{b'} = \overline{b}}>> Z @>{\overline{d'} =\overline{d}}>> X \langle1\rangle  @>{-\overline{a'}=-\overline{a}\langle1\rangle}>> Y\langle1\rangle
\end{CD}
$$
in $\uu$.

(TR3) Take a commutative diagram
$$
\begin{CD}
X @>{\overline{a}}>> Y @>{\overline{b}}>> Z @>{\overline{c}}>> X \langle1\rangle \\
@V{\overline{f}}VV @V{\overline{g}}VV @. @V{\overline{f} \langle1\rangle}VV \\
X' @>{\overline{a'}}>> Y' @>{\overline{b'}}>> Z' @>{\overline{c'}}>> X' \langle1\rangle
\end{CD}
$$
of triangles in $\uu$.
By the construction in \ref{subfactor3.2}, there exists a (not necessarily commutative) diagram
$$
\begin{CD}
X @>{a}>> Y @>{b}>> Z @>{d}>> X[1] \\
@V{f}VV @V{g}VV @. @V{f[1]}VV \\
X' @>{a'}>> Y' @>{b'}>> Z' @>{d'}>> X'[1]
\end{CD}
$$
of triangles in $\tt$.
Since $\overline{ag}=\overline{fa'}$ holds,
$ag-fa'$ factors through $\dd$. Since $a$ is $\dd$-monic,
there exists $s\in[\dd](Y,X')$ such that $ag-fa'=as$.
Replacing $g$ by $g-s$, we can assume that $ag=fa'$ holds.
By the axiom (TR3) for $\tt$, there exists $h\in\tt(Z,Z')$ which makes the diagram
$$
\begin{CD}
X @>{a}>> Y @>{b}>> Z @>{d}>> X[1] \\
@V{f}VV @V{g}VV @V{h}VV @V{f[1]}VV \\
X' @>{a'}>> Y' @>{b'}>> Z' @>{d'}>> X'[1]
\end{CD}
$$
commutative.
Thus the assertion follows from \ref{subfactor3.4}(3).

(TR4) Let
$$
X \stackrel{\overline{a}}{\longrightarrow}Y \stackrel{\overline{b}}{\longrightarrow}Z \stackrel{\overline{c}}{\longrightarrow} X\langle1\rangle, \quad  Y \stackrel{\overline{d}}{\longrightarrow} U \stackrel{\overline{e}}{\longrightarrow} V \stackrel{\overline{f}}{\longrightarrow} Y\langle1\rangle, \quad  X \stackrel{\overline{ad}}{\longrightarrow} U \stackrel{\overline{g}}{\longrightarrow} W \stackrel{\overline{h}}{\longrightarrow} X\langle1\rangle
$$
be triangles in $\uu$.
By \ref{subfactor3.4}(1),  we may assume that $a\in\tt(X,Y)$ and $d\in\tt(Y,U)$ are $\dd$-monic. Then $ad\in\tt(X,U)$ is also $\dd$-monic.
By \ref{subfactor3.4}(2), we have triangles
$$
X \stackrel{a}{\longrightarrow} Y \stackrel{b}{\longrightarrow} Z \stackrel{c'}{\longrightarrow}X[1], \quad Y \stackrel{d}{\longrightarrow} U \stackrel{e}{\longrightarrow} V \stackrel{f'}{\longrightarrow}Y[1], \quad  X \stackrel{ad}{\longrightarrow} U \stackrel{g}{\longrightarrow} W \stackrel{h'}{\longrightarrow} X[1]
$$
in $\tt$ which induce the above triangles in $\uu$.
By the octahedral axiom in $\tt$, we have commutative diagrams
$$
\begin{CD}
X @>{a}>> Y @>{b}>> Z @>{c'}>> X[1] \hspace{24pt} \\
@| @VdVV @VVV \hspace{-24pt} @|  \\
X @>{ad}>> U @>{g}>> W @>{h'}>> X[1] \hspace{24pt}  @. X @>{ad}>> U @>{g}>> W @>{h'}>> X[1] \\
@. @V{e}VV @ViVV @. @V{a}VV @| @ViVV @V{a[1]}VV \\
@. V @= V @. @. Y @>{d}>> U @>{e}>> V @>{f'}>> Y[1] \\
@. @V{f'}VV @VVV \\
@. Y[1] @>{b[1]}>> Z[1] \\
\end{CD}
$$
of triangles in $\tt$.
By \ref{subfactor3.4}(3), we have commutative diagrams of triangles:
$$
\begin{CD}
X @>{\overline{a}}>> Y @>{\overline{b}}>> Z @>{\overline{c}}>> X \langle1\rangle \hspace{24pt} \\
@| @V{\overline{d}}VV @VVV \hspace{-24pt} @|  \\
X @>{\overline{ad}}>> U @>{\overline{g}}>> W @>{\overline{h}}>> X \langle1\rangle \hspace{24pt}  @.
X @>{\overline{ad}}>> U @>{\overline{g}}>> W @>{\overline{h}}>> X \langle1\rangle \\
@. @V{\overline{e}}VV @V{\overline{i}}VV @. @V{\overline{a}}VV @| @V{\overline{i}}VV  @V{\overline{a}\langle1\rangle}VV \\
@. V @= V @. @. Y @>{\overline{d}}>> U @>{\overline{e}}>> V @>{\overline{f}}>> Y \langle1\rangle \\
@. @V{\overline{f}}VV @VVV \\
@. Y \langle1\rangle @>{\overline{b}\langle1\rangle}>> Z \langle1\rangle
\end{CD}
$$
\rightline{\qed}

We denote by $\overline{(-)}:\zz\to\uu$ the natural functor.
We make the following observation which we use later.

\begin{prop}\label{subfactor mutation}
\begin{itemize}
\item[(1)]
If $\cc$ is a 2-rigid subcategory of $\zz$, then so is $\overline{\cc}$ as a subcategory of $\uu$.

\item[(2)]
Assume that $(\xx,\yy)$ is a $\cc$-mutation pair in $\tt$ such that $\xx\vee\yy\subset\zz$ and $\dd\subset\cc$. Then $(\overline{\xx},\overline{\yy})$ is a $\overline{\cc}$-mutation pair in $\uu$.
\end{itemize}
\end{prop}

\begin{pf}
(1) Fix $X\in\cc$. Then the assertion follows from the exact sequence
$$
\begin{CD}
\tt(\cc,X) @>\cdot\alpha_X>> \tt(\cc,D_X) @>\cdot\beta_X>> \tt(\cc,X\langle1\rangle) @>\cdot\gamma_X>> \tt(\cc,X[1])=0.
\end{CD}
$$

(2) For any $X\in\xx$  (resp. $Y\in\yy$),  take a triangle $X \stackrel{f}{\to} C \stackrel{g}{\to} Y \to X[1]$ in $\tt$ with $Y\in\yy$ and a right $\cc$-approximation $g$ (resp. $X\in\xx$ and a left $\cc$-approximation $f$). Then $X \stackrel{\overline{f}}{\to} C \stackrel{\overline{g}}{\to} Y \to X\langle1\rangle$ is a triangle in $\uu$ with $Y\in\overline{\yy}$ and a right $\overline{\cc}$-approximation $\overline{g}$ (resp. $X\in\overline{\xx}$ and a left $\overline{\cc}$-approximation $\overline{f}$). Thus $(\overline{\xx},\overline{\yy})$ is a $\overline{\cc}$-mutation pair.
\qed
\end{pf}

We shall need the following lemma as well.

\begin{lemma}\label{subfactor3.5}
For any $X\in\zz$ and $i\ge0$, there exists a triangle
$$
\begin{CD}
C_i @>{\beta_X^{(i)}}>> X \langle i\rangle @>{\gamma_X^{(i)}}>> X[i] @>>> C_i[1] \\
\end{CD}
$$
in $\tt$ with $C_i\in\dd*\dd[1]*\cdots*\dd[i-1]$ and with $\beta_X^{(i)}$ being $\dd$-epic .
\end{lemma}

\begin{pf}
We shall show the assertion by induction on  $i$.
Assume that we have a triangle
$$
\begin{CD}
C_{i-1} @>{\beta_X^{(i-1)}}>> X\langle{i-1}\rangle @>{\gamma_X^{(i-1)}}>> X[i-1] @>>> C_{i-1}[1]
\end{CD}
$$
in $\tt$ with $C_{i-1}\in\dd*\dd[1]*\cdots*\dd[i-2]$.
Take a triangle
$$
\begin{CD}
X \langle {i-1}\rangle @>{\alpha_{X\langle i-1\rangle}}>> D_{X\langle i-1\rangle} @>>> X \langle i\rangle @>>>X \langle i-1\rangle[1].
\end{CD}
$$
By the octahedral axiom, we have the following commutative diagram.
{\small\
$$
\begin{CD}
@. X \langle i \rangle [-1] @= X \langle i\rangle[-1] \\
@. @VVV @VVV @. \\
C_{i-1} @>{\beta_X^{(i-1)}}>> X \langle i-1\rangle @>{\gamma_X^{(i-1)}}>> X[i-1] @>>> C_{i-1}[1] \\
@| @V{\alpha_{X\langle i-1\rangle}}VV @VVV @| \\
C_{i-1} @>>> D_{X\langle i-1\rangle} @>>> C_i @>>> C_{i-1}[1] \\
@. @VVV @V{\beta_X^{(i)}}VV @. \\
@. X \langle i\rangle @= X \langle i \rangle @. \\
\end{CD}
$$}
Then $\beta_X^{(i)}$ is $\dd$-epic, and we obtain a desired triangle $X\langle i\rangle[-1]\to X[i-1]\to C_i\to X\langle i\rangle$ in $\tt$.
\qed
\end{pf}

In the rest of this section we consider a special case of the setting above.
For this,  let $\tt$ be a triangulated category with a Serre functor $\sss$ and let $\zz\supset\dd$ be $\sss_n$-subcategories of $\tt$ with $n\ge2$.  In the rest we assume that one of the following conditions (A) or (B) holds.
\begin{itemize}
\item[(A)] $\dd$ is a functorially finite $n$-rigid subcategory of $\tt$ and
$$\displaystyle\zz = \bigcap_{i=1}^{n-1}\dd[-i]^\perp= \bigcap_{i=1}^{n-1}{}^\perp\dd[i]\supset\dd.$$

\item[(B)] $n=2$ and $\zz$ and $\dd$ satisfy the conditions (Z1) and (Z2).
\end{itemize}

\begin{prop}\label{subfactor4.1}
$\zz$ and $\dd$ satisfy the conditions (Z1) and (Z2).
\end{prop}

\begin{pf}
We only have to consider the case (A).
(Z1) is obvious. We will show (Z2).
For any $X\in\zz$  (resp. $Y\in\zz$),  take a triangle $X \stackrel{f}{\to} D \stackrel{g}{\to} Y \to X[1]$ with a left $\dd$-approximation $f$ (resp. right $\dd$-approximation $g$). Then one can easily check that $Y\in\zz$ and $g$ is a right $\dd$-approximation ($X\in\zz$ and $f$ is a left $\dd$-approximation). Thus $(\zz,\zz)$ is a $\dd$-mutation pair.
\qed
\end{pf}

Now the following theorem is the main result of this section. This kind of reduction from $\tt$ to $\uu$ was applied in \cite{KR} and \cite{CK} to cluster categories.

\begin{thm}\label{subfactor4.2}
Under the assumption (A) or (B), the subfactor category  $\uu = \zz/[\dd]$  forms a triangulated category with a Serre functor $\sss':=\sss_n \circ \langle n \rangle$.
In particular, if $\tt$ is $n$-Calabi-Yau, then so is $\uu$.
\end{thm}

To prove the theorem we need a lemma.

\begin{lemma}\label{subfactor4.3}
For any $i \ (0< i \le n)$, we have a functorial monomorphism
$$
\uu (Y, X \langle i\rangle) \stackrel{\cdot \gamma_X^{(i)}}{\longrightarrow} \tt(Y,X[i])
$$
for any $X,Y \in \zz$.
This is an isomorphism if  $0 < i < n$.
\end{lemma}

\begin{pf}
By \ref{subfactor3.5}, we have an exact sequence
$$
\begin{CD}
\tt (Y,C_i) @>{\cdot \beta_X^{(i)}}>> \tt (Y,X\langle i\rangle) @>{\cdot \gamma_X^{(i)}}>>  \tt (Y,X[i])  @>{\cdot c}>> \tt (Y,C_i[1]),
\end{CD}
$$
where  $\beta _X^{(i)}$ is $\dd$-epic, hence  $\Im (\cdot \beta_X^{(i)}) \supset [\dd] (Y,X \langle i\rangle)$.
Note that  $\tt (Y, C_i)=[\dd] (Y, C_i)$, since  $C_i \in \dd*\dd[1]*\cdots*\dd[i-1]$ and  $\tt ( \zz, \dd[1]*\cdots*\dd[i-1]) = 0$.
Consequently, we have $\Im ( \cdot \beta_X^{(i)})=[\dd] (Y, X \langle i\rangle)$.
Thus the first assertion follows.
If $i\neq n$,  then  $(Y, C_i[1]) = 0$ holds.
Thus $\uu (X,Y \langle i\rangle) \to \tt (X,Y[i])$ is an isomorphism.
\qed
\end{pf}

Now we proceed to the proof of Theorem \ref{subfactor4.2}.

\vspace{12pt}
\noindent{\bf Proof of \ref{subfactor4.2}}
By \ref{subfactor3.3} and \ref{subfactor4.1}, we have only to show that $\sss_n \circ \langle n\rangle$ is a Serre functor for  $\uu$.
Take any $X, Y \in \zz$.
By \ref{subfactor3.5}, we have a triangle
$C_n \stackrel{\beta_X^{(n)}}{\longrightarrow} X \langle n \rangle \stackrel{\gamma_X^{(n)}}{\longrightarrow} X [n] \stackrel{c}{\to} C_n[1]$
in $\tt$ with $C_n \in \dd * \dd[1] * \cdots * \dd[n-1]$.
We have a commutative diagram of exact sequences
$$
\begin{CD}
D\tt (Y, C_n[1]) @>{D(\cdot c)}>>  D\tt (Y, X[n]) @>{D(\cdot \gamma_X^{(n)})}>> D\tt (Y, X \langle n \rangle) @>{D(\cdot \beta_X^{(n)})}>> D\tt(Y,C_n) \\
@V{\wr}VV @V{\wr}VV @V{\wr}VV @V{\wr}VV  \\
\tt (C_n[1], \sss Y) @>{c\cdot}>> \tt (X[n], \sss Y) @>{\gamma_X^{(n)}\cdot}>>
\tt (X \langle n\rangle, \sss Y) @>{\beta_X^{(n)}\cdot}>> \tt(C_n,\sss Y),
\end{CD}
$$
which gives a functorial isomorphism
\begin{equation}\label{sf*}
\Im ( \gamma_X^{(n)} \cdot) \simeq D\Im (\cdot \gamma_X^{(n)}) = D\uu (Y, X\langle n\rangle),
\end{equation}
by \ref{subfactor4.3}.

Since $C_n[1] \in \dd[1]*\cdots*\dd[n]$ and $\tt( \dd[1]*\cdots*\dd[n-1], \sss Y) = 0$ holds by $\sss \zz = \zz[n]$,
we have $\tt (C_n[1], \sss Y) = [\dd[n]] (C_n[1], \sss Y)$.
Thus  $\Im (c \cdot ) \subset [\dd[n]] (X[n],\sss Y)$ holds.

On the other hand, since $\beta_X^{(n)}$ is $\dd$-epic by \ref{subfactor3.5},
$(\beta_X^{(n)} \cdot)=D(\cdot \beta_X^{(n)})$ is a monomorphism if $Y \in \dd$.
Thus $c$ is a $(\sss\dd)$-monic.
Since $\sss \dd = \dd[n]$ holds,  we have $\Im (c \cdot)=[\dd[n]] ( X[n], \sss Y)$ for any $X,Y\in\zz$.
Thus we have a functorial isomorphism
\begin{equation}\label{sf**}
\Im ( \gamma_X^{(n)} \cdot) = (\tt/[\dd[n]]) (X[n], \sss Y) \simeq \uu (X, \sss_nY).
\end{equation}
By  (\ref{sf*}) and (\ref{sf**}), we have proved the assertion.
\qed

\begin{thm}\label{subfactor4.4}
Under the assumption (A), the correspondence $\cc \mapsto \overline{\cc} := \cc/[\dd]$ gives
\begin{itemize}
\item[(1)] a one-one correspondence between $n$-cluster tilting subcategories of $\tt$ containing $\dd$  and $n$-cluster tilting subcategories of $\uu$, and

\item[(2)] a one-one correspondence between $n$-rigid $\sss_n$-subcategories of $\tt$ containing $\dd$  and  $n$-rigid $\sss'_n$-subcategories of $\uu$ for $\sss'_n:=\sss'\circ\langle -n\rangle$.
\end{itemize}
\end{thm}

\begin{pf}
Obviously, any $n$-rigid subcategory $\cc$ of $\tt$ containing $\dd$ is contained in $\displaystyle\zz =\bigcap_{i=1}^{n-1}\dd[-i]^\perp= \bigcap_{i=1}^{n-1}{}^\perp \dd[i]$.
Let $\cc$ be any subcategory of $\zz$ containing $\dd$.
By virtue of  \ref{subfactor4.3}, we have
\begin{eqnarray*}
\bigcap_{i=1}^{n-1} \overline{\cc \langle i\rangle}^\perp  = \overline {\bigcap_{i=1}^{n-1}\cc[i]^\perp}\ \ \ \mbox{and}\ \ \ \bigcap_{i=1}^{n-1}{}^\perp \overline{\cc \langle i\rangle}  =  \overline{\bigcap_{i=1}^{n-1}{}^\perp \cc[i]}.
\end{eqnarray*}
By \ref{subfactor4.2}, $\cc$ is an $\sss_n$-subcategory of $\tt$ if and only if $\overline{\cc}$ is an $\sss'_n$-subcategory of $\uu$.

Thus it is sufficient to verify the functorially finiteness. Since $\dd$ is a functorially finite subcategory of $\tt$, so is $\zz$ by \ref{torsion} and its dual.
Thus $\cc$ is a functorially finite subcategory of $\tt$ if and only if $\cc$ is so as a subcategory of $\zz$ if and only if $\overline{\cc}$ is so as a subcategory of $\uu=\zz/[\dd]$.
\qed
\end{pf}

In the next section we shall use the following simple observation which asserts that AR $(n+2)$-angles are preserved by the natural functor $\cc \to \overline{\cc}$.

\begin{prop} \label{subfactor4.5}
Under the assumption (A), let  $\cc$  be  an $n$-cluster tilting subcategory of $\tt$ containing $\dd$ and $\overline{\cc}:=\cc/[\dd]$.
If $X\in\cc$ does not have a direct summand in $\ind\dd$ and
$$
\sss_n X \stackrel{b_n}{\longrightarrow} C_{n-1} \ang{X_{n-1}} {f_{n-1}}{a_{n-1}}{b_{n-1}} C_{n-2} \ang{X_{n-2}}{f_{n-2}}{a_{n-2}}{b_{n-2}} \cdots \cdots
\ang{X_2}{f_2}{\ \ a_2}{b_2}C_1\ang{X_1}{f_1}{\ \ a_1}{b_1} C_0 \stackrel{a_0}{\longrightarrow} X
$$
is an AR $(n+2)$-angle of $X$ in $\cc$,
then its image
$$
\sss_nX \stackrel{\overline{b}_n}{\longrightarrow} C_{n-1} \ang{X_{n-1}}{\overline{f}_{n-1}}{\overline{a}_{n-1}}{\overline{b}_{n-1}} C_{n-2} \ang{X_{n-2}}{\overline{f}_{n-2}}{\overline{a}_{n-2}}{\overline{b}_{n-2}} \cdots \cdots \ang{X_2}{\overline{f}_2}{\ \ \overline{a}_2}{\overline{b}_2} C_1 \ang{X_1}{\overline{f}_1}{\ \ \overline{a}_1}{\overline{b}_1} C_0 \stackrel{\overline{a}_0}{\longrightarrow} X
$$
is an AR $(n+2)$-angle of $X$ in $\overline{\cc}$.
\end{prop}

\begin{pf}
Since each $a_i$ is a $\dd$-epic, each $X_i$ in fact belongs to $\zz$ by \ref{subfactor3.4}.
Obviously,  $\overline{a_0}$  is a sink map of $X$ in $\overline{\cc}$,  and $\overline{a_n}$ is a minimal right $\overline{\cc}$-approximation of $X_i$.
\qed
\end{pf}

\section{Mutation of $n$-cluster tilting subcategories}

Throughout this section, let $\tt$ be a triangulated category with Serre functor $\sss$.

\begin{thm}\label{2.1}
Let $\dd$ be a functorially finite $n$-rigid $\sss_n$-subcategory of $\tt$ and $(\xx,\yy)$ a $\dd$-mutation pair.
\begin{itemize}
\item[(1)] $\xx$ is an $n$-cluster tilting subcategory of $\tt$ if and only if so is $\yy$.
\item[(2)] $\xx$ is an $n$-rigid $\sss_n$-subcategory of $\tt$ if and only if so is $\yy$.
\end{itemize}
\end{thm}

\begin{pf}
Put $\displaystyle\zz:=\bigcap_{i=1}^{n-1}\dd[-i]^\perp=\bigcap_{i=1}^{n-1}{}^\perp\dd[i]$.
Then $\dd$ and $\zz$ satisfy the condition (A).
By \ref{subfactor4.2}, $\uu:=\zz/[\dd]$ is a triangulated category with the shift  functor $\langle1\rangle$ and the Serre functor $\sss'$.
By \ref{subfactor mutation}, $(\overline{\xx},\overline{\yy}):=(\xx/[\dd],\yy/[\dd])$ forms a $0$-mutation pair.
Thus $\overline{\yy}=\overline{\xx}\langle1\rangle$ holds.
In particular, $\overline{\xx}$ is an $n$-cluster tilting subcategory (resp. $n$-rigid $\sss'_n$-subcategory) of $\uu$ if and only if so is $\overline{\yy}$.
On the other hand, by \ref{subfactor4.4}, $\overline{\xx}$ (resp. $\overline{\yy}$) is an $n$-cluster tilting subcategory (resp. $n$-rigid $\sss'_n$-subcategory) of $\uu$ if and only if
$\xx$ (resp. $\yy$) is an $n$-cluster tilting subcategory (resp. $n$-rigid $\sss_n$-subcategory) of $\tt$. Thus the assertions follow.
\qed
\end{pf}

\begin{defn}
We call a functorially finite $n$-rigid $\sss_n$-subcategory $\dd$ of $\tt$ {\it almost complete $n$-cluster tilting} if there exists an $n$-cluster tilting subcategory $\cc$ of $\tt$ such that $\dd \subset \cc$ and $\ind\cc-\ind\dd$ consists of a single $\sss_n$-orbit. When $\tt$ is $n$-Calabi-Yau, then $\sss_n$ is the identity functor and these conditions are equivalent to  $\dd \subset \cc$ and $\#(\ind\cc-\ind\dd)=1$.

For an $\sss_n$-subcategory $\cc$ of $\tt$, we denote by $\sss_n\backslash(\ind\cc)$ the set of $\sss_n$-orbits in $\ind\cc$.
\end{defn}

For a given almost complete $n$-cluster tilting subcategory $\dd$ of $\tt$, we are interested in the set of $n$-cluster tilting subcategories containing $\dd$. Especially, it is natural to ask how many $n$-cluster tilting subcategories of $\tt$ contain $\dd$.
Partial results were given in \cite{Buan, Geiss, Thomas}.

For the case $n=2$, we can give the following complete answer.


\begin{thm}\label{2complements}
Any almost complete $2$-cluster tilting subcategory $\dd$ of $\tt$ is contained in exactly two $2$-cluster tilting subcategories $\xx$ and $\yy$ of $\tt$. Both $(\xx,\yy)$ and $(\yy,\xx)$ form $\dd$-mutation pairs.
\end{thm}

\begin{pf}
By \ref{2.1}, we only have to show the former assertion since $(\xx,\xx)$ never forms a $\dd$-mutation pair.
Put $\zz:=\dd[-1]^\perp={}^\perp\dd[1]$. Then $\uu:=\zz/[\dd]$ forms a triangulated category by \ref{subfactor3.3}.

Replacing $(\tt,\dd)$ by $(\uu,0)$, we can assume that $0$ is a almost complete $2$-cluster tilting subcategory of $\tt$ by \ref{subfactor4.4}. Thus we may assume that $\cc$ is a 2-cluster tilting subcategory of $\tt$ such that $\ind\cc$ consists of a single $\sss_2$-orbit. Then $\cc$ and $\cc[1]$ are distinct 2-cluster tilting subcategories of $\tt$ because $(\cc,\cc[1])=0$. We only have to show that any $X\in\ind\tt$ satisfying $(\sss_2^\ell X,X[1])=0$ for any $\ell\in\Z$ belongs to $\cc$ or $\cc[1]$.

By \ref{1.7}, there exists a triangle $C_1\stackrel{f}{\to}C_0\to X\to C_1[1]$ with $C_i\in\cc$ and $f\in J_{\tt}$. For any $\ell\in\Z$, chasing the commutative diagram
\[\begin{array}{ccccccc}
&&&&(\sss_2^\ell X[-1],X)=0&&\\
&&&&\uparrow&&\\
(\sss_2^\ell C_1,C_1)&\stackrel{\cdot f}{\to}&(\sss_2^\ell C_1,C_0)&\to&(\sss_2^\ell C_1,X)&&\\
&&\uparrow^{\sss_2^\ell f\cdot}&&\uparrow&&\\
&&(\sss_2^\ell C_0,C_0)&\to&(\sss_2^\ell C_0,X)&\to&(\sss_2^\ell C_0,C_1[1])=0,
\end{array}\]
we have that $(\sss_2^\ell C_1,C_1)\oplus(\sss_2^\ell C_0,C_0)\stackrel{{\cdot f\choose \sss_2^\ell f\cdot}}{\longrightarrow}(\sss_2^\ell C_1,C_0)\to0$ is exact. Since $f\in J_{\tt}$, we have $(\sss_2^\ell C_1,C_0)=J_{\tt}(\sss_2^\ell C_1,C_0)$ for any $\ell\in\Z$. This implies either $C_0=0$ or $C_1=0$ because $\ind\cc$ consists of a single $\sss_2$-orbit. Thus either $X=C_1[1]\in\cc[1]$ or $X=C_0\in\cc$ holds.
\qed\end{pf}

We note that the analog of \ref{2complements} for $n$-cluster tilting subcategories does not hold in general. In fact our results \ref{main2'} and \ref{theorem2 in section 10} give a $(2n+1)$-Calabi-Yau triangulated category which contains infinitely many $(2n+1)$-cluster tilting subcategories with one indecomposable object.

We shall give in \ref{ncomplements} a sufficient condition for $\tt$ so that any almost complete $n$-cluster tilting subcategory is contained in exactly $n$ $n$-cluster tilting subcategories. The condition will be given in terms of AR $(n+2)$-angles.

\begin{defn}\label{2.2}
Let $\cc$ be an $n$-cluster tilting subcategory of $\tt$ and let
\vspace{4pt}
$$
\sss_nX\stackrel{b_n}{\longrightarrow}C_{n-1}\ang{X_{n-1}}{f_{n-1}}{a_{n-1}}{b_{n-1}}C_{n-2}\ang{X_{n-2}}{f_{n-2}}{a_{n-2}}{b_{n-2}}\cdots\cdots\ang{X_2}{f_2}{\ \ a_2}{b_2}C_1\ang{X_1}{f_1}{\ \ a_1}{b_1}C_0\stackrel{a_0}{\longrightarrow}X
$$
be an AR $(n+2)$-angle of $X_0 = X \in \ind \cc$.
\begin{itemize}
\item[(1)]
 We say that $X$ {\it has no loops in $\cc$} if
$\displaystyle\sss_n^\ell X\notin\add\bigoplus_{i=0}^{n-1}C_i$
holds for any $\ell \in \Z$ and for the terms $C_i$ in the AR $(n+2)$-angle of $X$.
We say that $\cc$ has {\it no loops} if any $X\in\ind\cc$ has no loops in $\cc$.
\item[(2)]
 Assume that $X\in\ind\cc$ has no loops. For $i\in \Z/n\Z$, define a subcategory $\mu_{X}^{i}(\cc)$ of $\tt$ by
$$
\ind\mu_{X}^{i}(\cc)=(\ind\cc-\{\sss_n^\ell X\ |\ \ell\in\Z \})\coprod\{\sss_n^\ell X_i \ |\ \ell \in \Z \},
$$
where $X_i$ is the $i$-th term in the AR $(n+2)$-angle of $X$.

Now define a subcategory $\dd$ of $\tt$ by
$$
\ind\dd=\ind\cc-\{\sss_n^\ell X\ |\ \ell\in\Z \}.
$$
Since $X$ has no loops, $\dd$ is a functorially finite subcategory of $\cc$ and of $\tt$. Moreover, $(\mu^1_X(\cc),\cc)$ forms a $\dd$-mutation pair because $a_0$ is a right $\dd$-approximation and $b_1$ is a left $\dd$-approximation.
\end{itemize}
\end{defn}

For the case $n=2$, $\mu^1_X(\cc)$ was studied in \cite{Buan} and \cite{Geiss} as a `categorical realization' of Fomin-Zelevinsky mutation.





\begin{thm}\label{2.3}
Let $\cc$ be an $n$-cluster tilting subcategory of $\tt$. Assume that $X\in\ind\cc$ has no loops in $\cc$.

\begin{itemize}
\item[(1)]
 $\mu_{X}^{i}(\cc)$ is an $n$-cluster tilting subcategory of $\tt$ for any $i\in \Z/n\Z$.
\item[(2)]
 Under the notation of \ref{2.2},
$$
\sss_nX_i\stackrel{\sss_nb_i}{\longrightarrow}\sss_nC_{i-1}\stackrel{\sss_nf_{i-1}}{\longrightarrow}\cdots\cdots\stackrel{\sss_nf_1}{\longrightarrow}\sss_nC_{0}\stackrel{(\sss_na_0)b_n}{\longrightarrow}C_{n-1}\stackrel{f_{n-1}}{\longrightarrow}\cdots\cdots\stackrel{f_{i+1}}{\longrightarrow}C_i\stackrel{a_i}{\longrightarrow}X_i
$$
is an AR $(n+2)$-angle of $X_i$ in $\mu_{X}^{i}(\cc)$.
\item[(3)]
 $\mu_X^i (\cc) = \mu_{X_{i-1}}^1 \circ \cdots \circ \mu_{X_1}^1 \circ \mu_X^1(\cc)$
holds for any $i \in \Z/n\Z$.

\end{itemize}
\end{thm}

\begin{pf}
If we prove (2) for the case $i=1$, then the whole assertion follows by induction on $i$.
Since $(\mu^1_X(\cc),\cc)$ forms a $\dd$-mutation pair, (1) for $i=1$ follows from \ref{2.1}. Now we will prove (2) for $i=1$, namely
$$
\sss_nX_1\stackrel{\sss_nb_1}{\longrightarrow}\sss_nC_{0}\stackrel{(\sss_na_1)b_0}{\longrightarrow}C_{n-1}\stackrel{f_{n-1}}{\longrightarrow}\cdots\cdots\stackrel{f_3}{\longrightarrow}C_2\stackrel{f_2}{\longrightarrow}C_1\stackrel{a_1}{\longrightarrow}X_1
$$
is an AR $(n+2)$-angle of $X_1$ in $\mu_X^1(\cc)$.

Since $X$ has no loops in $\cc$, each term in this complex lies in $\mu_X^1(\cc)$.
By \ref{1.10}, we have only to show that $\sss_nb_1$ is a source map of $\sss_nX_1$ in $\mu_X^1(\cc)$, or equivalently, $b_1$ is a source map of $X_1$ in $\mu_X^1(\cc)$.
Since
$$
(C_0,-)\stackrel{b_1\cdot}{\longrightarrow}J_{\tt}(X_1,-)\to0
$$
is exact as functors on $\dd$,  it is enough to show that this is exact
after evaluating on the set $\{\sss_n^\ell X_1\ |\ \ell\in \Z \}$.
To show this, take any $f\in J_{\tt}(X_1,\sss_n^\ell X_1)$.
Since $cf\sss_n^\ell b_1\in (X[-1], \sss_n^\ell C_0)=0$,
there exist $g$ and $h$ which make the following diagram commutative.
$$
\begin{CD}
C_0[-1]  @>{a_0[-1]}>>  X[-1]  @>{c}>>  X_1  @>{b_1}>>  C_0  @>{a_0}>>  X  \\
@V{g[-1]}VV   @V{h[-1]}VV  @V{f}VV  @V{g}VV  @V{h}VV  \\
\sss_n^\ell C_0[-1]  @>{\sss_n^\ell a_0[-1]}>>  \sss_n^\ell X[-1]  @>{\sss_n^\ell c}>>  \sss_n^\ell X_1  @>{\sss_n^\ell b_1}>>
\sss_n^\ell C_0  @>{\sss_n^\ell a_0}>>  \sss_n^\ell X \\
\end{CD}
$$
If $h$ is an isomorphism, then so is $g$ by the right minimality of $a_0$.
Thus $f$ is also an isomorphism, a contradiction to $f\in J_{\tt}$.
Thus we have  $h\in J_{\tt}$.
Since $\sss_n^\ell a_0$ is a sink map of $\sss_n^\ell X$ in $\cc$, there exists
$s \in (X, \sss _n^\ell C_0)$ such that $h=s\sss_n^\ell a_0$.
Since $cf= h[-1]\sss_n^\ell c = s[-1]\sss_n^\ell (a_0[-1]c) =0$,
$f$ factors through $b_1$.
Thus we have shown that $b_1$ is a source map of $X_1$ in $\mu_X^1(\cc)$.
\qed
\end{pf}

From the viewpoint of tilting theory, it is natural to ask how much information on $\tt$ one can recover from its $n$-cluster tilting subcategory. We study certain special cases in \ref{4.1} and \ref{almost no loop}.

\begin{prop}\label{4.1}
Let $\cc$ be an $n$-cluster tilting subcategory of $\tt$.
\begin{itemize}
\item[(1)] $J_{\cc}=0$ and $\cc[n]=\cc$ hold if and only if any $X\in\cc$ has an AR $(n+2)$-angle
$$\sss_n X \stackrel{}{\longrightarrow} 0 \longrightarrow 0 \longrightarrow \cdots\cdots\longrightarrow 0 \longrightarrow 0 \stackrel{}{\longrightarrow}X.$$

\item[(2)] If the conditions in (1) hold, then $\displaystyle\ind\tt=\coprod_{i=0}^{n-1}\ind\cc[i]$.
\end{itemize}
\end{prop}

\begin{pf}
(1) Fix $X\in\cc$. The $(n+2)$-angle
$$
\sss_n X \stackrel{}{\longrightarrow} 0 \ang{X[1-n]}{}{}{} 0\longrightarrow \cdots\cdots\longrightarrow 0 \ang{X[-1]}{}{}{} 0 \stackrel{}{\longrightarrow}X
$$
is an AR $(n+2)$-angle if and only if $J_{\cc}(-,X)=0$ and $(\cc,X[-i])=0$ for $0<i<n$ if and only if $J_{\cc}(-,X)=0$ and $X[-n]\in\cc$. Thus (1) holds.

(2) By \ref{n-cluster torsion}, we have $\tt=\cc*\cc[1]*\cdots*\cc[n-1]$. Using \ref{direct summands} repeatedly, we have $\cc*\cdots*\cc[n-1]=\cc\vee\cdots\vee\cc[n-1]$ because $J_{\cc}=0$ and $(\cc[i],\cc)=0$ hold for $0<i<n$.
Since $\cc[i]\cap\cc[j]=0$ for any $0\le i\neq j<n$, we have the assertion.
\qed
\end{pf}



\begin{lemma}\label{4.2}
Let $\cc$ be an $n$-cluster tilting subcategory of $\tt$. Assume that $X \in \ind \cc$ has no loops, and define subcategories $\dd$ and $\zz$ of $\tt$ by
\begin{eqnarray*}
\ind \dd = \ind \cc-\{\sss_n^\ell X\ |\ \ell\in\Z\}\ \ \ \mbox{and}\ \ \ \zz=\bigcap_{i=1}^{n-1}\dd[-i]^\perp=\bigcap_{i=1}^{n-1}{}^\perp\dd[i].
\end{eqnarray*}
\begin{itemize}
\item[(1)] $\sss_n\backslash(\ind\zz)=\sss_n\backslash(\ind\dd)\coprod\{X, X_1,\cdots,X_{n-1}\}$
holds, where  $X_1, \cdots , X_{n-1}$ are the terms appearing in an AR $(n+2)$-angle
$$
\sss_nX\stackrel{}{\to}C_{n-1}\ang{X_{n-1}}{}{}{}C_{n-2}\ang{X_{n-2}}{}{}{}\cdots\cdots\ang{X_2}{}{}{}C_1\ang{X_1}{}{}{}C_0\stackrel{}{\to}X.
$$

\item[(2)]  $\ind \cc-\{\sss_n^\ell X\ |\ \ell\in\Z\}$ is contained in exactly $n$ $n$-cluster tilting subcategories $\mu_{X}^{i} (\cc)$ ($i \in \Z/n\Z$) of $\tt$.
\end{itemize}
\end{lemma}

\begin{pf}
We know that $\uu = \zz/[\dd]$ is a triangulated category with a Serre functor $\sss^\prime=\sss\circ\langle n\rangle$ by \ref{subfactor4.2}, and $\overline{\cc} = \cc/[\dd]$ is an $n$-cluster tilting subcategory of $\uu$ by \ref{subfactor4.4}.
By definition, $\ind\overline{\cc}$ consists of a single $\sss^\prime_n$-orbit.
Since $X$ has no loops in $\cc$,
$$
\sss_n X \stackrel{}{\longrightarrow} 0 \longrightarrow 0 \longrightarrow \cdots\cdots\longrightarrow 0 \longrightarrow 0 \stackrel{}{\longrightarrow}X
$$
is an AR $(n+2)$-angle of $X$ in $\overline{\cc}$ by \ref{subfactor4.5}. Thus $X_i=X\langle i\rangle$ holds for $0<i<n$. Applying \ref{4.1}, we have $\sss_n\backslash(\ind \uu)=\coprod_{i=0}^{n-1}\sss^\prime_n\backslash(\ind \overline{\cc}\langle i\rangle)=\{ X, X_1, \cdots, X_{n-1} \}$.
Thus (1) holds. We obtain (2) by \ref{subfactor4.4}.
\qed
\end{pf}

Consequently, we have the following main result of this section.

\begin{thm}\label{4.3}
Let $\cc$ be an $n$-cluster tilting subcategory of $\tt$. If $X \in \ind \cc$ has no loops in $\cc$, then $\ind \cc-\{\sss_n^\ell X\ |\ \ell\in\Z\}$ is contained in exactly $n$ $n$-cluster tilting subcategories $\mu_{X}^{i} (\cc)$ ($i \in \Z/n\Z$) of $\tt$.
\end{thm}


\vskip.5em
As an immediate consequence, we obtain the following result.

\begin{cor}\label{ncomplements}
Assume that any $n$-cluster tilting subcategory of $\tt$ has no loops. Then any almost complete $n$-cluster tilting subcategory of $\tt$ is contained in exactly $n$ $n$-cluster tilting subcategories of $\tt$.
\end{cor}

The following analogue of \ref{4.1} will be used later.

\begin{prop}\label{almost no loop}
Let $\cc$ be a $(2n+1)$-cluster tilting subcategory of $\tt$.
\begin{itemize}
\item[(1)] $J_{\tt}(\cc[i],\cc)=0$ holds for $0\le i<n$ if and only if any $X\in\cc$ has an AR $(2n+3)$-angle
$$\sss_n X \stackrel{}{\longrightarrow} \stackrel{2n}{0}\longrightarrow \cdots \longrightarrow \stackrel{n+1}{0} \longrightarrow \stackrel{n}{C_n} \longrightarrow \stackrel{n-1}{0}\longrightarrow \cdots \longrightarrow \stackrel{0}{0} \stackrel{}{\longrightarrow}X.$$
\item[(2)] If the conditions in (1) hold, then $\displaystyle\ind\tt=(\coprod_{i=0}^{n-1}\ind(\cc[i]*\cc[i+n+1]))\coprod\ind\cc[n]$.
\end{itemize}
\end{prop}

\begin{pf}
One can show the former assertion as in the proof of \ref{4.1}.
Fix $X\in\cc$. From the AR $(2n+3)$-angle, we obtain a triangle $\sss_nX[n]\to C_n\to X[-n]\to\sss_nX[n+1]$. In particular, $\cc\subset\cc[n]*\cc[2n+1]$ holds. This implies $(\cc[i],\cc)=0$ for $n<i<2n$.

Using \ref{direct summands} repeatedly, we have $\cc*\cdots*\cc[n]=\cc\vee\cdots\vee\cc[n]$ and $\cc[n+1]*\cdots*\cc[2n]=\cc[n+1]\vee\cdots\vee\cc[2n]$ because $J_{\tt}(\cc[i],\cc)=0$ hold for $0\le i<n$.
By \ref{n-cluster torsion}, we have
$$
\tt=(\cc*\cdots*\cc[n])*(\cc[n+1]*\cdots*\cc[2n])=(\cc\vee\cdots\vee\cc[n])*(\cc[n+1]\vee\cdots\vee\cc[2n]).
$$
Since $J_{\tt}(\cc[i],\cc)=0$ for any $0\le i<n$ and $n<i<2n$, we have
$$
\tt=(\cc*\cc[n+1])\vee\cdots\vee(\cc[n-1]*\cc[2n])\vee\cc[n]
$$
by applying \ref{direct summands} repeatedly.


Moreover, for any $0<i\neq j<n$, one can easily show that any morphism from $\cc[i]*\cc[i+n+1]$ to $\cc[j]*\cc[j+n+1]$ is in $[\cc[i+n+1]][\cc[j]]\subset J_{\tt}$. Thus $(\cc[i]*\cc[i+n+1])\cap(\cc[j]*\cc[j+n+1])=0$. Similarly, we have $(\cc[i]*\cc[i+n+1])\cap\cc[n]=0$. Thus we have the latter assertion.
\qed
\end{pf}

\section{Subfactor abelian categories}

Throughout this section, let $\tt$ be a triangulated category.

\begin{notation}
Let $\dd_i\subset\cc$ ($i=0,1$) be subcategories of $\tt$.

\begin{itemize}
\item[(1)] We denote by $\mod(\cc;\dd_0,\dd_1)$ the full subcategory of $\Mod\cc$ consisting of all $F$ which has a projective resolution $(-,D_1)\to(-,D_0)\to F\to0$ with $D_i\in\dd_i$.

\item[(2)] In this section, we shall exclusively consider functors of the form $\fff:\dd_1\to\Mod\dd_2$. We say that $\fff$ {\it preserves the 2-rigidity}
if $\Ext^1_{\Mod\dd_2}(\fff X,\fff Y)=0$  holds whenever  $\tt(X,Y[1])=0$ for $X,Y\in\dd_1$.
\end{itemize}
\end{notation}

\begin{prop}\label{equivalencelemma}
Let $\xx$, $\yy$ and $\zz$ be full subcategories of $\tt$ such that $\xx\vee\yy[-1]\subset\zz\subset{}^\perp\yy$. Then the following  (1)--(3) hold.
\begin{itemize}
\item[(1)] There exists an equivalence $\fff:(\xx*\yy)/[\yy]\to\mod(\zz;\xx,\yy[-1])$ defined by $\fff T:=\tt(-,T)$.

\item[(2)] $\fff$ preserves $2$-rigidity.

\item[(3)] If $\xx=\yy[-1]=\zz$, then $\fff$ gives an equivalence $(\xx*\yy)/[\yy]\to\mod\zz$. If $\xx={}^\perp\yy$, then $\fff$ gives a fully faithful functor $\tt/[\yy]\to\mod\zz$.
\end{itemize}
\end{prop}

\begin{pf}
(1) For any $T\in\xx*\yy$, there exists a triangle
\begin{eqnarray}\label{hajime}
\begin{CD}
Y[-1] @>f>> X @>g>> T @>h>> Y
\end{CD}
\end{eqnarray}
with $X\in\xx$ and $Y\in\yy$. Then we have an exact sequence
\begin{eqnarray}\label{ue}
(-,T[-1])\stackrel{\cdot h[-1]}{\longrightarrow}(-,Y[-1])\stackrel{\cdot f}{\longrightarrow}(-,X)\stackrel{\cdot g}{\longrightarrow}(-,T)\stackrel{\cdot h}{\longrightarrow}(-,Y)
\end{eqnarray}
on $\tt$. Restricting to $\zz\ (\subset{}^\perp\yy)$, we have a projective resolution
\begin{eqnarray}\label{tugi}
\begin{CD}
(-,Y[-1]) @>\cdot f>> (-,X) @>>>  \fff T @>>> 0
\end{CD}
\end{eqnarray}
of the $\zz$-module $\fff T$. Applying $\Hom_{\Mod\zz}(-,\fff T')$ and using Yoneda's Lemma, we have a commutative diagram
{\small\begin{eqnarray}\label{diagram}
\begin{array}{ccccccccc}
0&\stackrel{}{\longrightarrow}&(\fff T,\fff T')&\stackrel{(\cdot g)\cdot}{\longrightarrow}&((-,X),\fff T')&\stackrel{(\cdot f)\cdot}{\longrightarrow}&((-,Y[-1]),\fff T')&&\\
&&&&\wr\downarrow&&\wr\downarrow&&\\
(Y,T')&\stackrel{h\cdot}{\longrightarrow}&(T,T')&\stackrel{g\cdot}{\longrightarrow}&(X,T')&\stackrel{f\cdot}{\longrightarrow}&(Y[-1],T')&\stackrel{h[-1]\cdot}{\longrightarrow}&(T[-1],T')
\end{array}
\end{eqnarray}}
of exact sequences. Since the image of $(h\cdot)$ is $(\tt/[\yy])(T,T')$ by $\xx\subset{}^\perp\yy$, we have an isomorphism $\Hom_{\Mod\zz}(\fff T,\fff T')\simeq(\tt/[\yy])(T,T')$. Thus $\fff$ is fully faithful.

Obviously $\fff T$ belongs to $\mod(\zz;\xx,\yy[-1])$. To prove that $\fff$  is dense, we fix an arbitrary $F\in\mod(\zz;\xx,\yy[-1])$. Then we can take a projective resolution
\begin{eqnarray*}
\begin{CD}
(-,Y[-1]) @>\cdot f>> (-,X) @>>> F @>>> 0
\end{CD}
\end{eqnarray*}
with $Y\in\yy$ and $X\in\xx$. Take a triangle
$Y[-1] \stackrel{f}{\to} X \stackrel{g}{\to} T \stackrel{h}{\to} Y$
in $\tt$. Then $T\in\xx*\yy$. Comparing with (\ref{tugi}), we have $F\simeq \fff T$.

(2) If $(T,T'[1])=0$, then we have an exact sequence
$$
0\to(\fff T,\fff T')\to((-,X),\fff T')\to((-,Y[-1]),\fff T')\to0
$$
from (\ref{diagram}). One can easily check that this implies $\Ext^1_{\Mod\zz}(\fff T,\fff T')=0$.

(3) The former assertion follows from that $\mod(\zz;\xx,\yy[-1])=\mod\zz$. The latter one follows from \ref{wakamatsu}.
\qed\end{pf}

Immediately from this,  we have the iterative construction of $\tt/[\cc]$ that reminds us of the iterated module categories.

\begin{thm}\label{equivalence3}
Let $\cc$ be an $n$-cluster tilting subcategory of $\tt$. Define a chain $0=\tt_0\subset\tt_1\subset\cdots\subset\tt_{n}=\tt$ of subcategories of $\tt$ by
$$
\tt_\ell:=\cc*\cc[1]*\cdots*\cc[\ell-1]
$$
(see \ref{direct summands}, \ref{n-cluster torsion}). For any $0<\ell<n$, we have an equivalence
$$
\fff_\ell:\tt_{\ell+1}/[\cc[\ell]]\to\mod(\tt_{\ell};\tt_{\ell},{\cc[\ell-1]})
$$
defined by $\fff_\ell T:=\tt(-,T)$, which preserves the 2-rigidity.
\end{thm}

\begin{pf}
Putting $\xx=\zz:=\tt_{\ell}$ and $\yy:=\cc[\ell]$ and applying \ref{equivalencelemma}, we have the assertion.
\qed
\end{pf}

As one of the applications of this theorem, we have the following result.

\begin{cor}\label{equivalence}
Let $\cc$ be an $n$-cluster tilting subcategory of $\tt$. For any $0<\ell<n$, we have
$$
\bigcap_{0<i<n,\ i\neq \ell}\cc[-i]^\perp=\displaystyle\bigcap_{0<i<n,\ i\neq n-\ell}{}^\perp\cc[i],
$$
and we denote this by $\zz_\ell$. Then there exists an equivalence $\fff_\ell:\zz_\ell/[\cc]\to\mod\cc$ defined by $\fff_\ell X:=\tt(-,X[\ell])$, which preserves the $2$-rigidity. In particular, $\zz_1/[\cc]\simeq\cdots\simeq\zz_{n-1}/[\cc]$ is an abelian category.
\end{cor}

\begin{pf}
Since $\cc^\perp=\cc[1]*\cdots*\cc[n-1]={}^\perp\cc[n]$ holds by \ref{n-cluster torsion}, we have the desired equality. Notice that $\mod\cc$ forms an abelian category by \ref{coherent}.
Since $\zz_{1}=\cc[-1]*\cc$ holds by \ref{n-cluster torsion}, the assertions for $\ell=1$ follows from \ref{equivalence3}. For any $X\in\zz_{\ell+1}$ (resp. $Y\in\zz_\ell$), take a triangle
$X \stackrel{f}{\to} C \stackrel{g}{\to} Y \to X[1]$
with a left $\cc$-approximation $f$ (resp. right $\cc$-approximation $g$). One can easily check that $Y\in\zz_{\ell}$ and $g$ is a right $\cc$-approximation ($X\in\zz_{\ell+1}$ and $f$ is a left $\cc$-approximation).
Thus $(\zz_{\ell+1},\zz_\ell)$ is a $\cc$-mutation pair. By \ref{mutationequivalence}, we have an equivalence $\ggg_{\ell+1}:\zz_{\ell+1}/[\cc]\simeq\zz_{\ell}/[\cc]$ defined by $\ggg_{\ell+1}X:=Y$. Then $\fff_{\ell+1}=\fff_\ell\circ\ggg_{\ell+1}$ holds, and the assertions follow.
\qed
\end{pf}

Restricting ourselves to 2-cluster tilting subcategories, we have the following result as a corollary to \ref{equivalence3}.
We should note that the claims (1) and (3) in the corollary already appear in the papers by Buan, Marsh and Reiten \cite{BMR}, Keller and Reiten \cite{KR} and by Koenig and Zhu \cite{KZ}.

\begin{cor}\label{2-cluster}
Let $\cc$ be a $2$-cluster tilting subcategory of $\tt$.
\begin{itemize}
\item[(1)] There exists an equivalence $\fff:\tt/[\cc[1]]\to\mod\cc$ defined by $\fff T:=\tt(-,T)$, which preserves the 2-rigidity. Thus $\tt/[\cc[1]]$ is an abelian category.

\item[(2)] Assume that $\tt$ is 2-Calabi-Yau and ${\rm gl.dim}(\mod\cc)\le1$. If $T\in\tt$ has no direct summand in $\ind\cc[1]$, then $\tt(T,T[1])=0$ if and only if $\Ext^1_{\mod\cc}(\fff T,\fff T)=0$.

\item[(3)] Assume that $\tt$ has a Serre functor $\sss$.
\begin{itemize}
\item[(i)] $\fff$ induces an equivalence between $\sss\cc$ and injective objects in $\mod\cc$.

\item[(ii)] Any injective object in $\mod\cc$ has projective dimension at most one.

\item[(iii)]  If $T\in\tt$ has no direct summand in $\ind\cc[1]$, then $\fff(\sss T[-1])\simeq\tau(\fff T)$, where $\tau$ is the Auslander-Reiten translation (see \ref{coherent}).
\end{itemize}
\end{itemize}
\end{cor}

\begin{pf}
(1) Immediate from \ref{equivalence3} and \ref{coherent}.

(2) We have only to show the `if' part. We follow the notation in the proof of \ref{equivalencelemma}, where $\xx=\zz:=\cc$ and $\yy:=\cc[1]$ in this case. Assume $\Ext^1_{\mod\cc}(\fff T,\fff T)=0$. Since $T$ has no direct summand in $\ind\cc[1]$, we have $h\in J_{\tt}$ in the triangle (\ref{hajime}). Thus $\Im(\cdot h[-1])\subset J_{\tt}(-,Y[-1])$ holds in (\ref{ue}). Since ${\rm gl.dim}(\mod\cc)\le1$, the injection $\Im (\cdot h[-1])\to (-,Y[-1])$ splits. Thus we have $(\cdot h[-1])=0$ on $\cc$, and
$$
\begin{CD}
0 @>>> (-,Y[-1]) @>\cdot f>> (-,X) @>>> \fff T @>>> 0
\end{CD}
$$
is exact. Applying $\Hom_{\mod\cc}(-,\fff T)$ and considering the commutative diagram (\ref{diagram}), we have that $(X,T)\stackrel{f\cdot}{\to}(Y[-1],T)\to0$ is exact. Thus $0\to(T,Y[1])\stackrel{\cdot f[2]}{\longrightarrow}(T,X[2])$ is exact, and we have a commutative diagram
$$
\begin{array}{ccccccc}
(X,Y)=0&&&&&&\\
\uparrow^{\cdot h}&&&&&&\\
(X,T) &\stackrel{f\cdot}{\longrightarrow}& (Y[-1],T)&\stackrel{}{\longrightarrow}0&&&\\
\uparrow^{\cdot g}&&\uparrow^{\cdot g}&&&&\\
(X,X) &\stackrel{f\cdot}{\longrightarrow}& (Y[-1],X) &\stackrel{h[-1]\cdot}{\longrightarrow}& (T[-1],X) &\stackrel{g[-1]\cdot}{\longrightarrow}& (X[-1],X)=0\\
&&\uparrow^{\cdot f}&&\uparrow^{\cdot f}&&\\
&& (Y[-1],Y[-1]) &\stackrel{h[-1]\cdot}{\longrightarrow}& (T[-1],Y[-1]) &&
\end{array}
$$
of exact sequences. Chasing this diagram, we have an exact sequence $(T[-1],Y[-1])\stackrel{\cdot f}{\to}(T[-1],X)\to0$. Since $(T,Y)\stackrel{\cdot f[1]}{\to}(T,X[1])\stackrel{\cdot g[1]}{\to}(T,T[1])\stackrel{\cdot h[1]}{\longrightarrow}(T,Y[1])\stackrel{\cdot f[2]}{\longrightarrow}(T,X[2])$ is exact, we have $(T,T[1])=0$.

(3)(i)(ii) Fix any $X\in\cc$. Then (i) follows from $(-,\sss X)\simeq D(X,-)$.
By \ref{n-cluster torsion}, there exists a triangle
$\sss X[-1] \to C_1 \stackrel{f}{\to} C_0 \to \sss X$ with $C_i\in\cc$.
Now (ii) follows from an exact sequence
$$
\begin{CD}
0=D(X[-1],-) @>>> (-,C_1) @>\cdot f>> (-,C_0) @>>> D(X,-) @>>> (-,C_1[1])=0.
\end{CD}
$$

(iii) By \ref{n-cluster torsion}, there exists a triangle $C_1\stackrel{f}{\to} C_0\to T\to C_1[1]$ with $C_i\in\cc$. Then we have a minimal projective resolution $(-,C_1) \stackrel{\cdot f}{\to} (-,C_0) \to \fff T \to 0$ because $T$ has no direct summand in $\ind\cc[1]$. The assertion follows from the commutative diagram
$$
\begin{array}{rcccccc}
0=(-,\sss C_0[-1]) &\longrightarrow& \fff(\sss T[-1]) &\stackrel{}{\longrightarrow}& (-,\sss C_1) &\stackrel{\cdot\sss f}{\longrightarrow}& (-,\sss C_0)\\
&&&&\wr\downarrow&&\wr\downarrow\\
0&\longrightarrow& \tau(\fff T) &\stackrel{}{\longrightarrow}& D(C_1,-) &\stackrel{D(f\cdot)}{\longrightarrow}& D(C_0,-).\mbox{\qed}
\end{array}
$$
\end{pf}

\section{Application of Kac's theorem}


Throughout this section, let $k$ be an algebraically closed field.
We shall apply Kac's theorem in the representation theory of quivers to a classification of rigid objects in triangulated categories. Let us start with recalling basic facts concerning finite dimensional hereditary $k$-algebras. Such an algebra $H$ is Morita equivalent to a path algebra $kQ$ of certain quiver $Q$ without oriented cycles. We denote by $K_0(H)$ the Grothendieck group of $\mod H$. For $X\in\mod\Lambda$, we denote by $\underline{\dim}X$ the class of $X$ in $K_0(\mod H)$. Thus $K_0(H)$ is an abelian group with the basis $\underline{\dim}S_1,\cdots,\underline{\dim}S_n$, where $S_1,\cdots,S_n$ are simple $H$-modules.

We denote by $\Delta$ (resp. $\Delta_+$) the set of roots (resp. positive roots) associated to the underlying graph of the quiver $Q$. We denote by $\alpha_1,\cdots,\alpha_n$ the simple roots, and by $W$ the Weyl group. Note that $\alpha_i=\underline{\dim}S_i$, and we have a $W$-invariant quadratic form $q:K_0(H)\to\Z$, defined as  $q(\underline{\dim}X)=\dim_k\End_H(X)-\dim_k\Ext^1_H(X,X)$.
We call an element of
$$
\Delta^{\rm re}:=\bigcup_{i=1}^nW\alpha_i\ \ \ \ \ ({\rm resp.}\ \  \Delta^{\rm re}_+:=\Delta^{\rm re}\cap\Delta_+)
$$
a {\it real root} (resp. {\it positive real root}). Then a root $d$ is real if and only if $q(d)>0$ if and only if $q(d)=1$ (e.g. \cite[1.9]{Kac}). Let us recall the following result due to Kac \cite{Kac} (see also \cite{Gabriel-Roiter}).

\begin{thm}\label{Kac}
$\underline{\dim}$ gives a surjective map $\ind(\mod H)\to\Delta_+$. If $d\in\Delta^{\rm re}_+$, then there exists unique $X\in\ind(\mod H)$ such that $\underline{\dim}X=d$.
\end{thm}

Consequently, we have $\underline{\dim}X\in\Delta^{\rm re}_+$ for any rigid $X\in\ind(\mod H)$ because $q(\underline{\dim}X)=\dim_k\End_H(X)-\dim_k\Ext^1_H(X,X)=\dim_k\End_H(X)>0$. We call $d\in\Delta^{\rm re}_+$ a real {\it Schur} root if the corresponding $X\in\ind(\mod H)$ satisfies $\endm_H(X)=k$. We denote by $\Delta^{\rm reS}_+$ the set of real Schur roots. For Schur roots, we refer to \cite{Schofield, Derksen-Weyman} for example. Immediately, $\underline{\dim}$ gives a bijection from the set of indecomposable rigid $H$-modules to $\Delta^{\rm reS}_+$.

Applying this to the classification of rigid objects, we have the following result.

\begin{thm}\label{Kaccor}
Let $\tt$ be a triangulated category with a 2-cluster tilting subcategory $\cc=\add M$. Assume that $H:=\End_{\tt}(M)$ is hereditary, and we denote by $\Delta^{\rm reS}_+$ the set of real Schur roots associated to $H$. Then the map
$$
\begin{CD}
\tt @>\tt(M,-)>> \mod H @>\underline{\dim}>> K_0(H)
\end{CD}
$$
induces an injective map from 2-rigid objects in $\ind\tt-\ind\cc[1]$ to $\Delta^{\rm reS}_+$. This is bijective provided $\tt$ is 2-Calabi-Yau.
\end{thm}

\begin{pf}
Both assertions follow from \ref{2-cluster} and \ref{Kac}.
\qed
\end{pf}


%


\begin{ex}\label{kronecker}
We consider the quiver $Q:$\ \begin{picture}(32,5)
\put(0,3){\circle*{4}}
\put(0,3){\vector(1,0){28}}
\put(10,5){\scriptsize$m$}
\put(30,3){\circle*{4}}
\end{picture},
which we describe as having $m$ arrows.
In this case, we call $H:=kQ=\left(\begin{array}{cc}k&k^m\\ 0&k\end{array}\right)$ the {\it Kronecker algebra} of degree $m$. The corresponding Weyl group is generated by two elements
$s_1=\left(\begin{array}{cc}-1&0\\ m&1\end{array}\right)$ and
$s_2=\left(\begin{array}{cc}1&m\\ 0&-1\end{array}\right)$.
Thus $\Delta^{\rm re}_+$ consists of
$$
(1\ 0)s_2s_1s_2s_1\cdots s_{1\ \mbox{\scriptsize or }2}\ \ \ \mbox{ and }\ \ \ (0\ 1)s_1s_2s_1s_2\cdots s_{2\ \mbox{\scriptsize or }1}.
$$
There exist two indecomposable projective $H$-modules $P_0={k\choose 0}$ and $P_1={k^m\choose k}$, two indecomposable injective $H$-modules $I_0=D(0\ k)$ and $I_1=D(k\ k^m)$, and two simple $H$-modules $P_0$ and $I_0$.

Now we assume $m\ge2$.
It is easily shown \cite{Gabriel-Roiter} that the Auslander-Reiten quiver of $\mod H$ contains a preinjective component ${\bf C}_1$ and a preprojective component ${\bf C}_2$
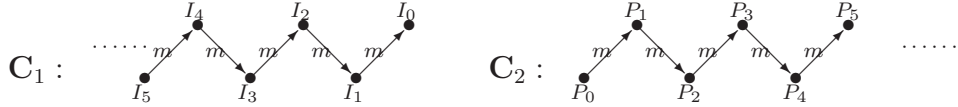
\begin{eqnarray*}
{\bf C}_1:\ \ \ \ \ \ \
\begin{picture}(110,20)
\put(-20,10){\scriptsize$\cdots\cdots$}
\put(0,0){\circle*{4}}
\put(-5,-8){\scriptsize$I_5$}
\put(0,0){\vector(1,1){18}}
\put(20,20){\circle*{4}}
\put(15,23){\scriptsize$I_4$}
\put(20,20){\vector(1,-1){18}}
\put(40,0){\circle*{4}}
\put(35,-8){\scriptsize$I_3$}
\put(40,0){\vector(1,1){18}}
\put(60,20){\circle*{4}}
\put(55,23){\scriptsize$I_2$}
\put(60,20){\vector(1,-1){18}}
\put(80,0){\circle*{4}}
\put(75,-8){\scriptsize$I_1$}
\put(80,0){\vector(1,1){18}}
\put(100,20){\circle*{4}}
\put(95,23){\scriptsize$I_0$}
\put(83,8){\scriptsize$m$}
\put(63,8){\scriptsize$m$}
\put(43,8){\scriptsize$m$}
\put(23,8){\scriptsize$m$}
\put(3,8){\scriptsize$m$}
\end{picture}&&
{\bf C}_2:\ \ \
\begin{picture}(140,25)
\put(0,0){\circle*{4}}
\put(-5,-8){\scriptsize$P_0$}
\put(0,0){\vector(1,1){18}}
\put(20,20){\circle*{4}}
\put(15,23){\scriptsize$P_1$}
\put(20,20){\vector(1,-1){18}}
\put(40,0){\circle*{4}}
\put(35,-8){\scriptsize$P_2$}
\put(40,0){\vector(1,1){18}}
\put(60,20){\circle*{4}}
\put(55,23){\scriptsize$P_3$}
\put(60,20){\vector(1,-1){18}}
\put(80,0){\circle*{4}}
\put(75,-8){\scriptsize$P_4$}
\put(80,0){\vector(1,1){18}}
\put(100,20){\circle*{4}}
\put(95,23){\scriptsize$P_5$}
\put(120,10){\scriptsize$\cdots\cdots$}

\put(83,8){\scriptsize$m$}
\put(63,8){\scriptsize$m$}
\put(43,8){\scriptsize$m$}
\put(23,8){\scriptsize$m$}
\put(3,8){\scriptsize$m$}
\end{picture}
\end{eqnarray*}
where all the modules are pairwise non-isomorphic and $\tau P_{i+2}=P_i$, $\tau I_i=I_{i+2}$ ($i\ge0$).
We have $\underline{\dim}P_0=(1,0)$, $\underline{\dim}P_1=(m,1)$, $\underline{\dim}I_0=(0,1)$ and $\underline{\dim}I_1=(1,m)$. Using Auslander-Reiten sequences, we have equalities$$
\underline{\dim}P_{i+1}-m(\underline{\dim}P_{i})+\underline{\dim}P_{i-1}=0,\ \ \ \ \
\underline{\dim}I_{i+1}-m(\underline{\dim}I_{i})+\underline{\dim}I_{i-1}=0
$$
for $i>0$. Therefore we obtain
\begin{eqnarray*}
\underline{\dim}P_{2i}=(1,0)(s_2s_1)^i,&\ \ \ &\underline{\dim}P_{2i+1}=(0,1)(s_1s_2)^is_1,\\
\underline{\dim}I_{2i}=(0,1)(s_1s_2)^i,&\ \ \ &\underline{\dim}I_{2i+1}=(1,0)(s_2s_1)^is_2.
\end{eqnarray*}
In particular, all elements in $\Delta^{\rm re}_+$ appear in ${\bf C}_1$ and ${\bf C}_2$. Consequently, any indecomposable rigid $H$-module is either $P_i$ or $I_i$ ($i\in\Z$).
\end{ex}

We apply \ref{kronecker} to give more explicit description of rigid objects.

\begin{cor}\label{2-cluster corollary}
Let $\tt$ be a $2$-Calabi-Yau triangulated category with a $2$-cluster tilting subcategory $\cc=\add(X\oplus Y)$. Assume $\End_{\tt}(X)=\End_{\tt}(Y)=k$, $(X,Y)=0$ and $\dim(Y,X)=m\ge2$.
\begin{itemize}
\item[(1)] $H:=\End_{\tt}(X\oplus Y)$ is a Kronecker algebra of degree $m$. We have an equivalence $\fff:=\tt(X\oplus Y,-):\tt/[\cc[1]]\to\mod H$ which preserves the 2-rigidity.

\item[(2)] The Auslander-Reiten quiver of $\tt$ contains a connected component
$$
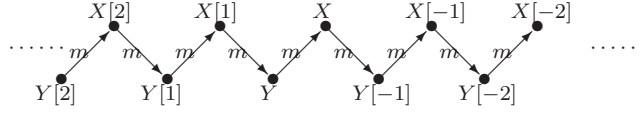
\begin{picture}(200,20)
\put(-20,10){\scriptsize$\cdots\cdots$}
\put(0,0){\circle*{4}}
\put(-10,-8){\scriptsize$Y[2]$}
\put(0,0){\vector(1,1){18}}
\put(20,20){\circle*{4}}
\put(10,23){\scriptsize$X[2]$}
\put(20,20){\vector(1,-1){18}}
\put(40,0){\circle*{4}}
\put(30,-8){\scriptsize$Y[1]$}
\put(40,0){\vector(1,1){18}}
\put(60,20){\circle*{4}}
\put(50,23){\scriptsize$X[1]$}
\put(60,20){\vector(1,-1){18}}
\put(80,0){\circle*{4}}
\put(75,-8){\scriptsize$Y$}
\put(80,0){\vector(1,1){18}}
\put(100,20){\circle*{4}}
\put(95,23){\scriptsize$X$}
\put(100,20){\vector(1,-1){18}}
\put(120,0){\circle*{4}}
\put(110,-8){\scriptsize$Y[-1]$}
\put(120,0){\vector(1,1){18}}
\put(140,20){\circle*{4}}
\put(130,23){\scriptsize$X[-1]$}
\put(140,20){\vector(1,-1){18}}
\put(160,0){\circle*{4}}
\put(150,-8){\scriptsize$Y[-2]$}
\put(160,0){\vector(1,1){18}}
\put(180,20){\circle*{4}}
\put(170,23){\scriptsize$X[-2]$}
\put(200,10){\scriptsize$\cdots\cdots$}

\put(163,8){\scriptsize$m$}
\put(143,8){\scriptsize$m$}
\put(123,8){\scriptsize$m$}
\put(103,8){\scriptsize$m$}
\put(83,8){\scriptsize$m$}
\put(63,8){\scriptsize$m$}
\put(43,8){\scriptsize$m$}
\put(23,8){\scriptsize$m$}
\put(3,8){\scriptsize$m$}
\end{picture}
$$
where all the modules are pairwise non-isomorphic.

\item[(3)] Any 2-rigid object in $\ind\tt$ is $X[i]$ or $Y[i]$ for some $i\in\Z$. Any 2-rigid object in $\tt$ is $X[i]^a\oplus Y[i]^b$ or $X[i]^a\oplus Y[i-1]^b$ for some $i\in\Z$ and $a,b\in\Z_{\ge0}$.

\item[(4)] Put $M_{2i}:=Y[-i]$, $M_{2i+1}:=X[-i]$ and $\cc_i:=\add(M_i\oplus M_{i+1})$ for $i\in\Z$. Then any 2-cluster tilting subcategory of $\tt$ is $\cc_i$ for some $i\in\Z$. They satisfy $\mu_{M_i}(\cc_i)=\cc_{i+1}$ and $\mu_{M_{i+1}}(\cc_i)=\cc_{i-1}$.
\end{itemize}
\end{cor}

\begin{pf}
(1) Immediate from our assumption and \ref{2-cluster}.

(3)(4) We use the notation in \ref{kronecker}. We have $\fff X=P_1$ and $\fff Y=P_0$. By \ref{2-cluster}, we have $\fff X[2]=I_0$, $\fff Y[2]=I_1$ and
\begin{eqnarray}\label{tau}
\fff X[-i]=\tau^{-i}P_1,\ \ \fff Y[-i]=\tau^{-i}P_0,\ \ \fff X[i+2]=\tau^iI_0,\ \ \fff Y[i+2]=\tau^iI_1
\end{eqnarray}
for $i\in\Z_{\ge0}$. Since any 2-rigid indecomposable $H$-module is either $P_i$ or $I_i$ by \ref{kronecker}, we have the former assertion of (3).

By \ref{2complements}, we have only to show that $\cc_i$ is a 2-cluster tilting subcategory of $\tt$. Since $\tt$ is 2-Calabi-Yau, our assumption implies that $\cc$ has AR 4-angles
$$
X \longrightarrow 0 \ang{X[1]}{}{}{} Y^m \longrightarrow X\ \ \ \mbox{ and } \ \ \ Y \longrightarrow X^m \ang{Y[-1]}{}{}{} 0 \longrightarrow Y.
$$
It follows from these 4-angles that $\cc_{-1}=\mu_{X}(\cc)$ and $\cc_1=\mu_{Y}(\cc)$, which are actually 2-cluster tilting subcategories of $\tt$ by virtue of \ref{2.3}. Since $\cc_{2i}=\cc[-i]$ and $\cc_{2i+1}=\cc_1[-i]$, we have that each $\cc_i$ is a 2-cluster tilting subcategory of $\tt$.
One can easily check $\mu_{M_i}(\cc_i)=\cc_{i+1}$ and $\mu_{M_{i+1}}(\cc_i)=\cc_{i-1}$.


(2) From the above 4-angles, $\tt$ has a triangle $X[1]\to Y^m\to X\to X[2]$. This is an AR triangle in $\tt$ because $\dim(X,X[2])=\dim(X,X)=1$. Similarly, $\tt$ has an AR triangle $Y\to X^m\to Y[-1]\to Y[1]$. Shifting these triangles, we have the assertion. All modules are pairwise non-isomorphic by (\ref{tau}).
\qed
\end{pf}

%

\begin{cor}\label{(2n+1)-cluster corollary}
Let $\tt$ be a $(2n+1)$-Calabi-Yau triangulated category with a $(2n+1)$-cluster tilting subcategory $\cc=\add X$ ($n>0$). Assume $\End_{\tt}(X)=k$, $(X[i],X)=0$ for $0<i<n$ and $\dim(X[n],X)=m\ge2$.
\begin{itemize}
\item[(1)] Put $\xx_{\ell}:=\cc[\ell]*\cc[\ell+n]*\cc[\ell+n+1]$ for $0\le \ell<n$. Then $\displaystyle\ind\tt=(\coprod_{\ell=0}^{n-1}(\ind\xx_\ell-\ind\cc[\ell+n+1]))\coprod\ind\cc[2n]$.

\item[(2)] $H:=\End_{\tt}(X\oplus X[n])$ is a Kronecker algebra of degree $m$. For $0\le \ell<n$, we have an equivalence $\fff_\ell:=\tt(X\oplus X[n],-[-\ell]):\xx_\ell/[\cc[\ell+n+1]]\to\mod H$ which preserves the 2-rigidity.

\item[(3)] The Auslander-Reiten quiver of $\tt$ contains a connected component
$$
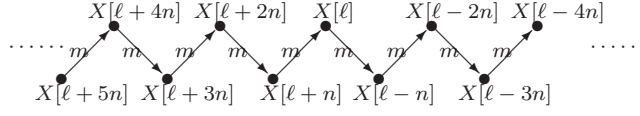
\begin{picture}(200,20)
\put(-20,10){\scriptsize$\cdots\cdots$}
\put(0,0){\circle*{4}}
\put(-10,-8){\scriptsize$X[\ell+5n]$}
\put(0,0){\vector(1,1){18}}
\put(20,20){\circle*{4}}
\put(10,23){\scriptsize$X[\ell+4n]$}
\put(20,20){\vector(1,-1){18}}
\put(40,0){\circle*{4}}
\put(30,-8){\scriptsize$X[\ell+3n]$}
\put(40,0){\vector(1,1){18}}
\put(60,20){\circle*{4}}
\put(50,23){\scriptsize$X[\ell+2n]$}
\put(60,20){\vector(1,-1){18}}
\put(80,0){\circle*{4}}
\put(75,-8){\scriptsize$X[\ell+n]$}
\put(80,0){\vector(1,1){18}}
\put(100,20){\circle*{4}}
\put(95,23){\scriptsize$X[\ell]$}
\put(100,20){\vector(1,-1){18}}
\put(120,0){\circle*{4}}
\put(110,-8){\scriptsize$X[\ell-n]$}
\put(120,0){\vector(1,1){18}}
\put(140,20){\circle*{4}}
\put(130,23){\scriptsize$X[\ell-2n]$}
\put(140,20){\vector(1,-1){18}}
\put(160,0){\circle*{4}}
\put(150,-8){\scriptsize$X[\ell-3n]$}
\put(160,0){\vector(1,1){18}}
\put(180,20){\circle*{4}}
\put(170,23){\scriptsize$X[\ell-4n]$}
\put(200,10){\scriptsize$\cdots\cdots$}

\put(163,8){\scriptsize$m$}
\put(143,8){\scriptsize$m$}
\put(123,8){\scriptsize$m$}
\put(103,8){\scriptsize$m$}
\put(83,8){\scriptsize$m$}
\put(63,8){\scriptsize$m$}
\put(43,8){\scriptsize$m$}
\put(23,8){\scriptsize$m$}
\put(3,8){\scriptsize$m$}
\end{picture}
$$
for $0\le\ell<n$, where all $X[i]$ ($i\in\Z$) are pairwise non-isomorphic.

\item[(4)] Any 2-rigid object in $\ind\tt$ is $X[i]^a$ for some $i\in\Z$ and $a\in\Z_{\ge0}$.

\item[(5)] Any $(2n+1)$-cluster tilting subcategory of $\tt$ is $\add X[i]$ for some $i\in\Z$.
\end{itemize}
\end{cor}

\begin{pf}
(1) Since $J_{\tt}(\cc[\ell+n],\cc[\ell+1])=0$ and $J_{\tt}(\cc[\ell+n+1],\cc[\ell+n+1])=0$, we have $\xx_\ell=(\cc[\ell]\vee\cc[\ell+n])*\cc[\ell+n+1]=(\cc[\ell]*\cc[\ell+n+1])\vee\cc[\ell+n]$ by \ref{direct summands}. Thus the assertion follows from \ref{almost no loop}.

(2) The former assertion follows from our assumption. We only have to show the latter assertion for $\ell=0$. Put $\xx=\zz:=\cc*\cc[n]=\cc\vee\cc[n]$ and $\yy:=\cc[n+1]$. Then $\xx_0=\xx*\yy$ and we have an equivalence
$$
\fff_0=\tt(X\oplus X[n],-):\xx_0/[\cc[n+1]]\to\mod(\cc\vee\cc[n];\cc\vee\cc[n],\cc[n])
$$
by \ref{equivalencelemma}. Since $J_{\tt}(\cc,\cc\vee\cc[n])=0$, we have $\mod(\cc\vee\cc[n];\cc\vee\cc[n],\cc[n])=\mod(\cc\vee\cc[n])\simeq\mod H$.

(3) We have only to show the assertion for $\ell=0$.
Since $X\in\cc$ has an AR $(2n+3)$-angle
$$\begin{CD}
X \stackrel{}{\longrightarrow} \stackrel{2n}{0} \longrightarrow \cdots \longrightarrow \stackrel{n+1}{0} \longrightarrow \stackrel{n}{X^m} \longrightarrow \stackrel{n-1}{0}\longrightarrow \cdots \longrightarrow \stackrel{0}{0} \stackrel{}{\longrightarrow}X\end{CD}
$$
by \ref{almost no loop}, there exists a triangle $X[n]\to X^m\to X[-n]\to X[n+1]$ in $\tt$. This is an AR triangle because $\dim(X[-n],X[n+1])=\dim(X,X)=1$. Shifting this, we have the component. We will show that all $X[i]$ ($i\in\Z$) are pairwise non-isomorphic in the proof of (4).

(4) We use the notation in \ref{kronecker}. More strongly, we will show that any 2-rigid object in $\ind\xx_\ell$ is $X[\ell-in]$ ($i\ge-1$) or $X[\ell+in+1]$ ($i\ge1$), and they satisfy
\begin{eqnarray}\label{F}
\begin{array}{cc}
\fff_\ell(X[\ell-2in])=\tau^{-i}P_1,&\fff_\ell(X[\ell-(2i-1)n])=\tau^{-i}P_0,\\
\fff_\ell(X[\ell+(2i+2)n+1])=\tau^iI_0,&\fff_\ell(X[\ell+(2i+3)n+1])=\tau^iI_1.
\end{array}
\end{eqnarray}

We may assume $\ell=0$. By a similar argument as in the proof of \ref{2-cluster}, we have
\begin{eqnarray}\label{2n-AR}
\fff_0(T[2n])=\tau(\fff_0 T)
\end{eqnarray}
if both $T$ and $T[2n]$ are contained in $\xx_0$ and $T$ has no direct summand $X[n+1]$.

First, we shall show that $X[-in]$ ($i\ge-1$) and $X[in+1]$ ($i\ge1$) in fact belong to $\xx_0$. We have $\cc[-n]\subset\cc*\cc[n+1]$ by the triangle in (3). Thus we have $\xx_0[-n]=(\cc[-n]*\cc[1])\vee\cc\subset(\cc*\cc[n+1]*\cc[1])\vee\cc=\cc*\cc[n+1]*\cc[1]$. Since $J_{\tt}(\cc[1],\cc[1]*\cc[n+2])=0$, we have $\cc*\cc[n+1]*\cc[1]=(\cc*\cc[n+1])\vee\cc[1]$ by \ref{direct summands}. Thus we have
$$
\xx_0[-n]\subset\xx_0\vee\cc[1].
$$
In particular, if $X[-in]\in\xx_0$ and $X[-in]\neq X[n+1]$, then $X[-(i+1)n]\in\xx_0$.

Assume $X[-jn]\in\xx_0$ for $-1\le j\le i$ and $X[-jn]\neq X[n+1]$ for $-1\le j<i$. If we show $X[-in]\neq X[n+1]$, then $X[-(i+1)n]\in\xx_0$ by the above observation, and we will have $X[-in]\in\xx_0$ for any $i\ge-1$ inductively. Thus we assume $X[-in]=X[n+1]$. Then $i>0$ holds by the $(2n+1)$-rigidity of $X$, and we have $X[(2-i)n]\in\xx_0$.
Thus we have $\fff_0(X[(2-i)n])=\fff_0(X[3n+1])=(X\oplus X[n],X[3n+1])=D(X[n],X\oplus X[n])=I_1$.
Using (\ref{2n-AR}) repeatedly, we have $P_1=\fff_0X=\tau^{i-2}\fff_0(X[(2-i)n])=\tau^{i-2}I_1$. This is a contradiction because $P_1$ and $I_1$ do not belong to one component.

We have $\cc[2n+1]\subset\cc*\cc[n+1]$ by the triangle in (3). By a similar argument as above, we have
$$
\xx_0[n]\subset\xx_0\vee\cc[2n].
$$
This implies that if $X[in+1]\in\xx_0$ and $X[in+1]\neq X[n]$, then $X[(i+1)n+1]\in\xx_0$. Since we have already shown that $X[in+1]\neq X[n]$ for $i\ge1$, we have $X[in+1]\in\xx_0$ for any $i\ge1$.

Now we prove the desired assertion. Since we have $\fff_0 X=P_1$, $\fff_0 X[n]=P_0$, $\fff_0 X[2n+1]=I_0$ and $\fff_0 X[3n+1]=I_1$, we obtain (\ref{F}). Since any 2-rigid indecomposable $H$-module is either $P_i$ or $I_i$ by \ref{kronecker} and $\fff_0$ preserves the 2-rigidity, any 2-rigid object in $\ind\xx_\ell$ is $X[\ell-in]$ ($i\ge-1$) or $X[\ell+in+1]$ ($i\ge1$).

(5) Immediate from (4).
\qed
\end{pf}


\section{Gorenstein quotient singularities}

We shall apply the arguments used in the preceding sections to Cohen-Macaulay modules over a quotient singularity.

For this, in the rest of this paper, let  $k$  be an algebraically closed field of characteristic zero.
And let  $G$  be a finite subgroup of $\GL (d, k)$ which is a small subgroup, i.e. $G$ does not contain any pseudo-reflections except the identity.

Let  $V$  be an $d$-dimensional $k$-vector space with basis $\{ x_1, \ldots , x_d\}$.
The group $G$  naturally acts on $V$  and the action can be extended to the action on the regular local ring  $S = k[[ x_1, \ldots ,x_d]]$, the completion of the symmetric algebra over  $V$.
We denote by  $R = S^G$ the invariant subring of  $S$  by this action of  $G$.

We assume the following assumptions on the group $G$:

\begin{itemize}
\item[(G1)]
 $G$  is a subgroup of  $\SL (d, k)$.
\item[(G2)]
 Any element $\sigma \not= 1$ of $G$ does not have eigenvalue 1.
\end{itemize}

It is known by \cite{Watanabe} that the condition (G1) is equivalent to that
$R$  is a Gorenstein ring.
It is also known that (G2) is a necessary and sufficient condition for  $R$
 to have at most an isolated singularity.
For the lack of references, first we note this fact in the following.

For each element  $\sigma \in G$ we let
$$
W_{\sigma} = (\sigma - 1)(V),
$$
which is a subspace of  $V$.
By definition,  $\sigma$  is a pseudo-reflection if and only if $\Kdim _k W_{\sigma} \leq 1$.
Therefore, since $G$ is small, we have  $\Kdim _k W_{\sigma} \geqq 2$  for $\sigma \not= 1$.
Using this notation, the following proposition holds.

\begin{prop}\label{eigen}
The singular locus of  $R$  is a closed subset of  $\Spec (R)$ defined by the ideal
$$
I = \bigcap _{\sigma \not= 1 \in G} (W_{\sigma} S \cap R).
$$
\end{prop}

As an immediate consequence of this proposition we have the following corollary.

\begin{cor}\label{coreigen}
The following three conditions are equivalent.
\begin{itemize}
\item[(1)]
$R$  is an isolated singularity.
\item[(2)]
$W_{\sigma} = V$  for any  $\sigma \in G$ with  $\sigma \not=1$.
\item[(3)]
Any element $\sigma \not= 1$ of $G$ does not have eigenvalue 1.
\end{itemize}
\end{cor}

Even though there is no good reference for the proof of Proposition \ref{eigen}, this seems to be a well-known fact.
We shall give an outline of the proof below, and the precise proof of each step is left to the reader.

\begin{pf}
Let  $P \in \Spec (S)$  and  $\gp = P \cap R \in \Spec (R)$.
We let
$$
T_P = \{ \sigma \in G \ | \ x^{\sigma} - x \in P \ \text{for} \ \forall x \in S\} \subseteq
H_P = \{ \sigma \in G \ | \ P^{\sigma} = P\} \subseteq G,
$$
which are respectively called the inertia group and the decomposition group for $P$.
As a first step, we can show that

\vspace{6pt}
\noindent
(1) $T_P = \{1 \}$  if and only if  $P \not\supseteq \bigcap _{\sigma \not= 1 \in G} W_{\sigma} S$.
\vspace{6pt}
\par\noindent
Note from (1) that  $G$  is small if and only if  $T_P =\{1\}$  for any prime ideal $P$ of height one.

Let  $P = P_1, P_2, \ldots , P_r$  be all the prime ideals of $S$ lying over $\gp$.
Note that $r = [G:H_P]$.
Since $R_{\gp} \subseteq S_{\gp}$  is a finite extension, taking the $\gp$-adic completion of this, we have a finite extension
$$
\widehat{R_{\gp}} \subseteq \widehat{S_{\gp}} = \widehat{S_{P}} \times \widehat{S_{P_2}} \times \cdots  \times \widehat{S_{P_r}}.
$$
The action of $G$  can be naturally extended to  $S_{\gp}$, hence to  $\widehat{S_{\gp}}$.
And we can easily see that  $\widehat{R_{\gp}} = \widehat{S_{\gp}}^G$.
As a second step of the proof, one can show that

\vspace{6pt}
\noindent
(2) $H_P$ naturally acts on  $\widehat{S_{P}}$  and  $(\widehat{S_{P}})^{H_P} = R_{\gp}$.
\vspace{6pt}

Let  $\kappa (\gp)$  (resp. $\kappa (P)$) be the residue field of the local ring  $\widehat{R_{\gp}}$ (resp. $\widehat{S_{P}}$).
By (2) one can show that the field extension  $\kappa (\gp) \subseteq \kappa (P)$  is a Galois extension with Galois group  $H_P/T_P$.
Using this observation we can prove the following.

\vspace{6pt}
\noindent
(3) If $T_P = \{ 1\}$, then  $R_{\gp}$  is a regular local ring.
\vspace{6pt}

In fact, if  $T_P =\{1\}$, then $[\kappa (P) : \kappa (\gp)] = |H_P| = \rank _{\widehat{R_{\gp}}} \widehat{S_P}$  by (2).
This shows that  $\widehat{R_{\gp}} \subseteq \widehat{S_P}$  is a flat extension, and since  $\widehat{S_P}$  is regular, it follows that  $\widehat{R_{\gp}}$  is also regular.

On the other hand, one can show the following by using (2).

\vspace{6pt}
\noindent
(4) If  $\widehat{R_{\gp}} \subseteq \widehat{S_P}$  is a flat extension, then  the equality $\length _{S_P} ( S_P/ \gp S_P) = |T_P|$ holds.
\vspace{6pt}

In fact, if  $\widehat{R_{\gp}} \subseteq \widehat{S_P}$  is flat, then
it is actually free and $\length _{S_P} ( S_P/ \gp S_P) = \rank _{\widehat{R_{\gp}}} \widehat{S_P} \ /\  [\kappa (P) : \kappa (\gp)]$, which shows (4).

Finally we can prove the following.

\vspace{6pt}
\noindent
(5) If  $G$  is a small subgroup and  if $R_{\gp}$  is a regular local ring,  then $T_P =\{ 1\}$.
\vspace{6pt}

In fact, if  $\height \ P =1$,  then  $\widehat{R_{\gp}} \subseteq \widehat{S_P}$  is flat,  and we have from (4) that  $\length _{S_P} ( S_P/ \gp S_P) = |T_P|=1$, since  $G$ is small.
This means that the extension  $\widehat{R_{\gp}} \subseteq \widehat{S_P}$ is unramified if $\height \ P =1$.
Then, by the purity of branch locus we see that  $\widehat{R_{\gp}} \subseteq \widehat{S_P}$ is unramified for a given $P$, since this is a flat extension.
It follows from (4) that $|T_P| = 1$.

Finally, summing all up, we have the following equivalences, and the proof is completed.
$$
R_{\gp} \ \text{is regular}
 \ \Leftrightarrow \
T_P = \{1\}
 \ \Leftrightarrow \
P \not\supseteq \bigcap _{\sigma \not= 1 \in G} W_{\sigma}S
 \ \Leftrightarrow \
\gp \not\supseteq I \qquad \text{\qed}
$$
\end{pf}

In the rest we assume that $G$ satisfies the conditions (G1) and  (G2).
Then it is equivalent to that $R$  is a Gorenstein local ring of dimension $d$  with an isolated singularity.
Under such a situation we are interested in the category of maximal Cohen-Macaulay modules  $\CM (R)$  over $R$ and its stable category  $\pCMR$.


\vspace{6pt}
For the present, we discuss some generalities about these categories.
For this purpose let  $A$  be a general Gorenstein complete local ring of dimension $d$.
We denote by  $\pCMA$  the stable category of  $\CM (A)$.
By definition, $\pCMA$ is the factor category  $\CM (A)/[A]$.
We denote the set of morphisms in  $\CM (A) /[A]$  by  $\pHom _A (M, N)$  for  any  $M, N \in \pCMA$.
Since  $A$  is a complete local ring,  note that  $M$  is isomorphic to  $N$ in  $\pCMA$  if and only if  $M \oplus P \cong N \oplus Q$  in  $\CM (A)$  for some projective $A$-modules  $P$  and  $Q$.

For any $A$-module  $M$,  we denote the first syzygy module of  $M$  by  $\Omega _A M$.
  We should note that  $\Omega _A M$  is uniquely determined up to isomorphism  as an object in  the stable category.
 The $n$th syzygy module  $\Omega _A ^n M$  is  defined inductively  by  $\Omega _A ^n M = \Omega _A (\Omega _A ^{n-1}M)$,  for any nonnegative integer $n$.
Since  $A$  is a Gorenstein ring, it is easy to see that the syzygy functor  $\Omega _A : \pCMA  \to  \pCMA$  is an autofunctor.
Hence, in particular, we can define the cosyzygy functor  $\Omega _A^{-1}$  on  $\pCMA$  which is the inverse of  $\Omega_A$.
We note from   \cite[2.6]{Happel}  that  $\pCMA$  is a triangulated category with shift functor  $[1] = \Omega _A ^{-1}$.

Now we remark one of the fundamental dualities called the {\it Auslander-Reiten duality}, which was essentially given by Auslander \cite[I.8.8, III.1.8]{Auslander2}.

\begin{thm}\label{11}
Let  $A$  be a Gorenstein local ring of dimension $d$ as above.
And suppose that  $A$  has only an isolated singularity.
Then, for any  $X, Y \in \pCMA$,   we have a functorial isomorphism
$$
\Ext ^d_A (\pHom _A (X, Y), A) \cong \pHom _A (Y, X[d-1]).
$$
Therefore the triangulated category  $\pCMA$  is $(d-1)$-Calabi-Yau.
\end{thm}

\begin{pf}
This follows from the isomorphism given in \cite[3.10]{Yoshino1}
\qed
\end{pf}


Now we return to the original setting, that is,  $R = S^G$  where the subgroup $G$ of  $\GL (d, k)$  satisfies the conditions  (G1) and  (G2).
We always assume $d \geq 3$  in the following.
As we have remarked above as a general  result,  $\pCMR$ is a triangulated category with shift functor  $[1] = \Omega _R ^{-1}$ and it is  $(d-1)$-Calabi-Yau.

For the following result we refer to \cite[2.5]{Iyama1} and \cite[6.1]{Iyama2}, where $(d-1)$-cluster tilting subcategories are called maximal $(d-2)$-orthogonal subcategories.

\begin{thm}\label{maxortho}
Under the above circumstances, $\add S$ is a $(d-1)$-cluster tilting subcategory of  $\pCMR$.
Moreover, the Koszul complex
$$
\begin{CD}
S @>>> S^{{d\choose d-1}} @>>> S^{{d\choose d-2}} @>>> \cdots @>>> S^{{d\choose 2}} @>>> S^{{d\choose 1}} @>>> S
\end{CD}
$$
of $S$ gives an AR $(d+1)$-angle in $\pCMR$.
 \end{thm}


Thanks to this theorem, we can apply results in previous sections.
For example, we immediately obtain the following corollary by virtue of \ref{2-cluster}.


\begin{cor}\label{corequiv3}
Assume $d = 3$ and put $\Lambda:=\pEnd_R(S)$. We have an equivalence $\Ext^1 _R(S,-):\CMR /[S]\to\mod\Lambda$, which preserves the 2-rigidity.
\end{cor}

\begin{pf}
Since $\add S$ is equivalent to the category of finitely generated projective $\Lambda$-modules, $\mod(\add S)$ is equivalent to $\mod\Lambda$.
\qed
\end{pf}

In particular, if all 2-rigid $\Lambda$-modules are known, then all 2-rigid Cohen-Macaulay $R$-modules are also known.

\vskip1em
In the rest of this section, we give an application of this equivalence.

Let $k$ be an algebraically closed field again.
In the representation theory of finite dimensional $k$-algebras, the concept of representation type plays an important role.
A famous dichotomy theorem of Drozd asserts that any finite dimensional $k$-algebra is
either of {\it tame representation type} or {\it wild representation type}, and not both (see \cite{Drozd, Crawley-Boevey}).
In short, the representation theory of algebras of tame representation type can be approximated by that of
a polynomial ring $k[x]$ in one variable, while the representation theory of algebras of wild representation type is as complicated as that of a polynomial ring $k\langle x,y\rangle$ in two variables.
It is known to be `hopeless' to classify indecomposable modules over algebras of wild representation type (see \cite[4.4.3]{hopeless} and references there).
A path algebra $kQ$ of a quiver $Q$ is of tame representation type if $Q$ is either a Dynkin or extended Dynkin diagram, and is of wild representation type otherwise.
We call $Q$ a {\it wild quiver} if it is neither Dynkin nor extended Dynkin. The notion of tameness and wildness
together with the dichotomy theorem was introduced for one dimensional reduced Cohen-Macaulay $k$-algebras by Drozd-Greuel \cite{Drozd-Greuel}.
Unfortunately the general definition of tameness and wildness is still not established for arbitrary dimension.

The following result asserts that the category $\CMR$ is more complicated than $\mod(kQ)$ for the path algebra $kQ$ of a wild quiver $Q$.
Especially it follows that $R=S^G$ is never of finite representation type if $G\neq1$. This contains a result due to Auslander and Reiten \cite{AR2} as a special case.

\begin{prop}\label{wild}
Let $G\neq1$ be a finite subgroup of $\SL(3,k)$. Assume that $R:=S^G$ is an isolated singularity. Then there exists a wild quiver $Q$ and a dense functor $\CMR\to\mod(kQ)$.
\end{prop}

\begin{pf}
We have an equivalence $\CMR/[S]\to\mod\Lambda$ by \ref{corequiv3}. Thus we only have to show that there exists a surjection $\Lambda\to kQ$ for a wild quiver $Q$ because we have a dense functor $(kQ)\otimes_{\Lambda}(-):\mod\Lambda\to\mod(kQ)$. It is shown in the next lemma.
\qed
\end{pf}

\begin{lemma}
Under the notation of \ref{wild}, there exists an idempotent $e$ of $\Lambda$ such that $\Lambda/\Lambda e\Lambda$ is isomorphic to a path algebra $kQ$ of a wild quiver $Q$.
\end{lemma}

\begin{pf}
Since $R=S^G$ is an isolated singularity, $G$ is a cyclic group by \cite{cyclic}. Up to conjugation,
we can put $G=\langle{\rm diag}(\zeta^a,\zeta^b,\zeta^c)\rangle$ for a primitive $m$-th root of unity $\zeta$ and integers $a,b,c$ satisfying
$a+b+c\in m\Z$ and $(a,m)=(b,m)=(c,m)=1$ by (G2). Moreover, without loss of generality, we can assume $a=1\le b\le c<m$.

Put $\Gamma:=S*G$. Then $\Gamma$ has a presentation as follows: The set of vertices is $\Z/m\Z$.  For each vertex $i\in\Z/m\Z$, we draw three arrows $x:i\to i+a$, $y:i\to i+b$ and $z:i\to i+c$. The relations are $xy=yx$, $yz=zy$ and $zx=xz$.

Let $e_i$ be the primitive idempotent of $\Gamma$ corresponding to the vertex $i\in\Z/m\Z$. Choose $e_0$ in the way such that $\Lambda=\Gamma/\Gamma e_0\Gamma$. Put $E:=\{1,2,\cdots,b,b+1,c+1\}\subset\Z/m\Z$. Define an idempotent $e$ of $\Lambda$ by $e:=\sum_{i\in(\Z/m\Z)\backslash E}e_i$. The presentation of $\Lambda/\Lambda e\Lambda=\Gamma/\Gamma e\Gamma$ is given by simply deleting the vertex $(\Z/m\Z)\backslash E$ from that of $\Gamma$.
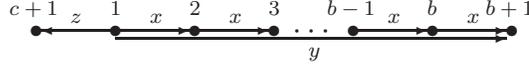
$$
\begin{picture}(180,10)
\put(0,0){\circle*{4}}
\put(-10,6){\scriptsize$c+1$}
\put(13,3){\scriptsize$z$}
\put(30,0){\circle*{4}}
\put(28,6){\scriptsize$1$}
\put(43,3){\scriptsize$x$}
\put(30,0){\vector(-1,0){28}}
\put(30,0){\vector(1,0){28}}
\put(30,-3){\vector(1,0){148}}
\put(60,0){\circle*{4}}
\put(58,6){\scriptsize$2$}
\put(73,3){\scriptsize$x$}
\put(60,0){\vector(1,0){28}}
\put(90,0){\circle*{4}}
\put(88,6){\scriptsize$3$}
\put(103,-10){\scriptsize$y$}
\put(97,-3){$\cdots$}
\put(120,0){\circle*{4}}
\put(110,6){\scriptsize$b-1$}
\put(133,3){\scriptsize$x$}
\put(120,0){\vector(1,0){28}}
\put(150,0){\circle*{4}}
\put(148,6){\scriptsize$b$}
\put(163,3){\scriptsize$x$}
\put(150,0){\vector(1,0){28}}
\put(180,0){\circle*{4}}
\put(170,6){\scriptsize$b+1$}
\end{picture}
$$
There are no relations among these arrows. Consequently, $\Lambda/\Lambda e\Lambda$ is the path algebra of a wild quiver. Thus the assertion follows.
\qed
\end{pf}


\section{Proof of Theorems \ref{main1} and \ref{main2}}

In the first half of this section we give a proof of Theorem \ref{main1} stated in the introduction.

For this, let  $G$  be a cyclic subgroup of $\GL (3, k)$  that is generated by  $\sigma = {\rm diag}(\omega, \omega, \omega)$, where  $\omega$ is a primitive cubic root of unity.
In this case, $G$  acts naturally on the ring  $S = k[[x, y, z]]$, and
the invariant subring   $R = S^{G}$  is the completion of the Veronese subring of dimension three and degree three:
$$
R = k [[ \{ \text{monomials of degree three in }  x, y, z \} ]]
$$
The action of  $G$  gives the $\ZC$-graded structure on $S$  in such a way that  $S = S_0 \oplus S_1 \oplus S_2$,  where each  $S_j$  is an $R$-module of semi-invariants that is defined as
$$
S_j = \{ f \in S \ | \ f^{\sigma} = \omega ^j f \}.
$$
Actually we have that $S_0=R$ and
$$
S_1 = (x,y,z)R, \quad  S_2 = (x^2, xy , y^2, yz, z^2 , zx)R.
$$

We shall prove Theorem \ref{main1} in the following form.

\begin{thm}\label{main1'}
Under the above circumstances, we set  $M_{2i}=\Omega_R^iS_1$  and  $M_{2i+1}=\Omega_R^iS_2$ for any $i\in\Z$. Then all $M_i$ ($i\in\Z$) are pairwise non-isomorphic.
\begin{itemize}
\item[(1)] There exists an equivalence $\CMR/[S]\to\mod\left(\begin{array}{cc}k&k^3\\ 0&k\end{array}\right)$.

\item[(2)] A maximal Cohen-Macaulay $R$-module is rigid if and only if it is isomorphic to $R^a \oplus M_i^b \oplus M_{i+1}^c$ for some $i\in\Z$ and $a, b, c\in\Z_{\ge0}$.

\item[(3)] Any 2-cluster tilting subcategory of $\pCMR$ is one of $\add(M_i\oplus M_{i+1})$ for $i\in\Z$.
\end{itemize}
\end{thm}


It is known that the $S_j \ (j \in \ZC)$  are in $\CMR$,  and in particular they are reflexive $R$-modules of rank one, whose classes form the divisor class group of  $R$.
As a result, $\{ S_j\}$ are all of the maximal Cohen-Macaulay modules of rank one over  $R$.


We need some results of computations for later use.

\begin{lemma}\label{4}
We have the following isomorphisms.
\begin{itemize}
\item[$(1)$]
$\Hom _R(S_i , S_j)  \cong  S_{j-i}$  for any  $i, j \in \ZC$.

\item[$(2)$]
$
\pHom _R (S_i, S_j) \cong
\begin{cases}
\  k  &(i=j=1 \ \ \text{or} \ \ i=j=2) \\
\  k^3  &(i=1, \ \ j=2) \\
\  0  &(i=2, \ \ j=1)\\
\end{cases}
$
\\
where we can take the multiplication maps  $S_1 \to S_2$  by $x, y, z$  as the $k$-basis of $\pHom _R(S_1, S_2)$.
\end{itemize}
\end{lemma}

\begin{pf}
(1) is obvious, and we omit its proof.
To show (2), let us take a Koszul complex of the sequence $\{ x, y, z\}$ in  $S$, that is an exact sequence of $\ZC$-graded $S$-modules.
$$
\begin{CD}
0 @>>> S @>>>S(-2)^3 @>>> S(-1)^3 @>f>> S @>>> k @>>> 0,
\end{CD}
$$
where  $f$  is a map given by the matrix
${}^t(x, y, z)$.
Taking the degree $1$ and $2$ part of this, we have exact sequences of $R$-modules.
$$
\begin{CD}
(i) \qquad  0 @>>> S_1 @>>>S_2 ^3 @>>> R^3 @>{f_1}>> S_1 @>>> 0, \\
(ii) \qquad  0 @>>> S_2 @>>> R ^3 @>>> S_1^3 @>{f_2}>> S_2 @>>> 0. \\
\end{CD}
$$
Note that $f_1$  gives a free cover of the $R$-module  $S_1$.
For any $j \in \ZC$, applying  $\Hom _R(S_j , \ \ )$ to this and using the isomorphisms of  $(1)$, we have the exact sequence
$$
\begin{CD}
S_{-j}^3 @>{f_{1*}}>> S_{-j+1}  @>>> \pHom _R(S_j, S_1) @>>> 0, \\
\end{CD}
$$
where $f_{1*}$  is the map induced by  $f$.
Hence  we have an isomorphism   $\pHom _R (S_j, S_1) \cong S_{-j+1}/(x, y, z)S_{-j}$.
  If  $j=1$, then we have  $\pHom _R (S_1, S_1) \cong k$, since  $(x, y ,z)S_2$  is the maximal ideal of $R =S_0$.
On the other hand, if  $j=2$, then we have
$\pHom _R(S_2, S_1) = 0$, since  $S_2 = (x, y, z)S_1$.

To show the rest of the lemma, we take a free presentation of  $S/(x, y, z)^2S$  as a $\ZC$-graded $S$-module.
$$
\begin{CD}
S(-2)^6  @>g>> S @>>> S/(x, y, z)^2S @>>> 0,
\end{CD}
$$
where  $g$  is a map given by the matrix
${}^t(x^2, xy, y^2, yz, z^2, zy)$.
Taking the degree $2$ part, we have a free presentation of  $S_2$  as an $R$-module as
$$
\begin{CD}
R^6 @>{g_2}>> S_2 @>>> 0, \\
\end{CD}
$$
For any $j \in \ZC$, applying  $\Hom _R(S_j , \ \ )$ again and using the isomorphisms of  $(1)$, we have an exact sequence
$$
\begin{CD}
S_{-j}^6 @>{g_{2*}}>> S_{-j+2}  @>>> \pHom _R(S_j, S_2) @>>> 0, \\
\end{CD}
$$
where $g_{2*}$  is the map induced by  $g$.
Hence  we have an isomorphism   $\pHom _R (S_j, S_2) \cong S_{-j+2}/(x, y, z)^2 S_{-j}$.
If  $j=1$, then we have  $\pHom _R (S_1, S_2) \cong S_1/(x, y, z)^2S_2 = S_1/(x, y, z)S_1 \cong k^3$.
On the other hand, if  $j=2$, then we have
$\pHom _R(S_2, S_2) = k$, since  $(x, y, z)^2 S_1$  is the maximal ideal of  $R = S_0$.
\qed\end{pf}

\vskip1em
\noindent
{\bf Proof of Theorem \ref{main1'} }
By \ref{maxortho}, $\cc:=\add S$ is a 2-cluster tilting subcategory of $\tt:=\pCMR$. Moreover, we have $\ind\cc=\{S_1,S_2\}$, $\End_{\tt}(S_1)=\End_{\tt}(S_2)=k$, $\tt(S_2,S_1)=0$ and $\dim\tt(S_1,S_2)=3$ by \ref{4}. Thus we obtain the assertions from \ref{2-cluster corollary} by putting $X:=S_2$ and $Y:=S_1$.
\qed
\vspace{12pt}


Now we proceed to the proof of  Theorem \ref{main2}.
For this purpose, we assume in the rest of the paper that
$G$  is a subgroup of  $\GL (4, k)$ generated by  ${\rm diag}(-1,-1,-1,-1)$,
and that  $G$  acts on the ring  $S = k[[x,y,z,w]]$ and  $R = S^G$.
As in the previous case, we decompose $S$  as the sum of modules of semi-invariants $S = S_0 \oplus S_1$, where  $S_0 = R$  and  $S_1 = (x, y, z, w)R$.
We repeat Theorem \ref{main2} in the introduction to give its proof.

\begin{thm}\label{main2'}
Under the above circumstances, all $\Omega^i_RS_1$ ($i\in\Z$) are pairwise non-isomorphic.
\begin{itemize}
\item[(1)] There exists an equivalence $\CMR/[S]\to\mod\left(\begin{array}{cc}k&k^6\\ 0&k\end{array}\right)$.

\item[(2)] A maximal Cohen-Macaulay $R$-module is rigid if and only if it is isomorphic to $R^a\oplus(\Omega^i_RS_1)^b$ for some $i\in\Z$ and $a,b\in\Z_{\ge0}$.

\item[(3)] Any 3-cluster tilting subcategory of $\pCMR$ is one of $\add\Omega^i_RS_1$ for $i\in\Z$.
\end{itemize}
\end{thm}

Note that  $\pCMR$  is a $3$-Calabi-Yau triangulated category in this case.
We need the following computational results for the proof.

\begin{lemma}\label{5}
Under the same circumstances, we have the following isomorphisms of  $k$-modules;
 $\underline{\End}_R(S_1)\simeq k$ and $\Ext^4_R(S_1,S_1)\simeq k^6$.
\end{lemma}

\begin{pf}
Similarly to the proof of \ref{4}(2),
we take a Koszul complex of the sequence $\{ x, y, z, w\}$ in  $S$, that is an exact sequence of $\Z / 2\Z$-graded $S$-modules.
$$
\begin{CD}
0 @>>> S @>>>S(-1)^4 @>>> S^6 @>>> S(-1)^4 @>f>> S @>>> k @>>> 0,
\end{CD}
$$
where  $f$  is a map given by the matrix
${}^t(x, y, z, w)$.
Taking the degree $1$ part of this, we have an exact sequence of $R$-modules.
$$
\begin{CD}
0 @>>> S_1 @>>> R^4 @>>> S_1 ^6 @>>> R^4 @>{f_1}>> S_1 @>>> 0 \\
\end{CD}
$$
Note that $f_1$  gives a minimal free cover of the $R$-module  $S_1$.
Hence we have  $\pHom _R(S_1, S_1) \simeq R/(x, y, z, w)R \simeq k$  as in the proof of \ref{4}(2).

On the other hand, from the above sequence we have a short exact sequence
$$
\begin{CD}
0 @>>> \Omega _R ^{-1} S_1 @>>> S_1 ^6 @>>> \Omega _R S_1 @>>> 0, \\
\end{CD}
$$
which is actually an Auslander-Reiten sequence, since
$$
\Ext _R^1(\Omega _RS_1, \Omega _R ^{-1}S_1) \simeq \Ext _R ^3(S_1, S_1) \simeq D\pHom _R(S_1, S_1) \simeq  k.
$$
It follows that  $\Ext ^4_R(S_1, S_1) \simeq D\pHom _R (S_1, \Omega _R S_1) \simeq D\pHom _R (S_1, S_1^6) \simeq k ^6$.
\qed
\end{pf}

\vspace{6pt}
\noindent
{\bf Proof of Theorem \ref{main2'}}
By \ref{maxortho}, $\cc:=\add S$ is a 3-cluster tilting subcategory of $\tt:=\pCMR$. Moreover, we have $\ind\cc=\{S_1\}$ and $\End_{\tt}(S_1)=k$ and $\dim\tt(S_1[1],S_1)=\dim\tt(S_1,S_1[4])=6$ by \ref{5}. Thus we obtain the assertions from \ref{(2n+1)-cluster corollary} by putting $n:=1$ and $X:=S_2$.
\qed


\section{Non-commutative examples}

In this section, we construct examples of $(2n+1)$-Calabi-Yau triangulated categories where
all $(2n+1)$-cluster tilting subcategories are known.


Let $R$ be a complete local Gorenstein isolated singularity of dimension $d$ and $M\in\CMR$ a generator such that $\add M$ is a $(d-1)$-cluster tilting subcategory of $\pCMR$.
Take a decomposition $M=M_1\oplus M_2$ as an $R$-module such that $M_1$ is a generator. Put
$$\Lambda:=\endm_R(M_1)\ \ \mbox{ and }\ \ N:=\Hom_R(M_1,M_2).$$
We denote by $\CML$ the category of $\Lambda$-modules which are Cohen-Macaulay as $R$-modules, and by $\pCML:=\CML/[\Lambda]$ the stable category.
We have the following analog of \ref{11} and \ref{maxortho}.

\begin{thm}\label{non-commutative}
Under the above circumstances, $\pCML$ is a $(d-1)$-Calabi-Yau triangulated category with a
$(d-1)$-cluster tilting subcategory $\add N$.
\end{thm}

\begin{pf}
Put $\Gamma:=\endm_R(M)$.
By \cite[5.2.1]{Iyama2}, we have $\Gamma\in \CMR$ and ${\rm gl.dim} \Gamma=d$.
Since $\Lambda$ is a direct summand of $\endm_R(M)$ as an $R$-module,
we have $\Lambda\in\CML$.
By \cite[2.4]{IyamaReiten}, $\Lambda$ is isomorphic to $\Hom_R(\Lambda,R)$ as an $(\Lambda,\Lambda)$-module. Thus $\pCML$ forms a $(d-1)$-Calabi-Yau triangulated category by
\cite[2.6]{Happel} and \cite[I.8.8, III.1.8]{Auslander2}.
Since $M_1$ is a generator, the functor $\Hom_R(M_1,-):\mod R\to\mod\Lambda$ is fully faithful.
Putting $N':=\Hom_R(M_1,M)$, we have $\endm_\Lambda(N')\simeq\Gamma$.
Again by \cite[5.2.1]{Iyama2}, we have that
$\add N=\add N'$ is a $(d-1)$-cluster tilting subcategory of $\pCML$.
\qed
\end{pf}

We apply the above construction to some special case.

Let $k$ be an algebraically closed field of characteristic zero and
$G$ a cyclic subgroup of $\GL(2n+2,k)$ with $n\ge0$ generated by $\sigma=
{\rm diag}(\zeta,\cdots,\zeta)$, where $\zeta$ is a primitive $(n+1)$-st root of unity.
Then the assumptions (G1) and (G2) in section 8 are satisfied.
Put $S:=k[[x_1,\cdots,x_{2n+2}]]$ and $R:=S^G$.
We observed in section 8 that $R$ is a Gorenstein isolated singularity and
$\pCMR$ forms a $(2n+1)$-Calabi-Yau triangulated category
with a $(2n+1)$-cluster tilting subcategory $\add S$.
Putting $S_j:=\{f\in S\ |\ f^\sigma=\zeta^jf\}$, we have a decomposition
$\displaystyle S=\bigoplus_{j=0}^nS_j$ as an $R$-module.
Put $\displaystyle M_1:=\bigoplus_{j=0}^{n-1}S_j$ and $M_2:=S_n$ to get $\Lambda:=\endm_R(M_1)$ and $N:=\Hom_R(M_1,M_2)$.
By \ref{non-commutative}, $\pCML$ is a $(2n+1)$-Calabi-Yau triangulated category with a
$(2n+1)$-cluster tilting subcategory $\add N$.

\begin{thm}\label{theorem2 in section 10}
Under the above circumstances, all $\Omega^i_\Lambda N$ ($i\in\Z$) are pairwise non-isomorphic.
\begin{itemize}
\item[(1)] A maximal Cohen-Macaulay $\Lambda$-module is rigid if and only if it is isomorphic to $P\oplus(\Omega^i_\Lambda N)^a$ for some projective $\Lambda$-module $P$, $i\in\Z$ and $a\in\Z_{\ge0}$.

\item[(2)] Any $(2n+1)$-cluster tilting subcategory of $\pCML$ is one of $\add\Omega^i_\Lambda N$ for $i\in\Z$.
\end{itemize}
\end{thm}

\begin{pf}
By \ref{maxortho}, we have an AR $(2n+3)$-angle
$$
\begin{CD}
S_n @>>> S_0^{{2n+2\choose 2n+1}} @>>> S_1^{{2n+2\choose 2n}} @>>> \cdots @>>> S_{n-2}^{{2n+2\choose 2}} @>>> S_{n-1}^{{2n+2\choose 1}} @>>> S_n
\end{CD}
$$
in the $(2n+1)$-cluster tilting subcategory $\add S$ of $\pCMR$.
Applying $\Hom_R(M_1,-)$, we have an AR $(2n+3)$-angle
$$
\begin{CD}
N @>>> 0 @>>> \cdots @>>> 0 @>>> N^{{2n+2\choose n+1}} @>>> 0 @>>> \cdots @>>> 0 @>>> N
\end{CD}
$$
in the $(2n+1)$-cluster tilting subcategory $\add N$ of $\pCML$.
By \ref{almost no loop} and \ref{(2n+1)-cluster corollary}, we have the assertion.
\qed
\end{pf}


\end{document}